\documentclass{amsart}
\usepackage{amsmath}
\usepackage{amsthm}
\usepackage{amssymb}
\usepackage{graphics}

\newcommand{\R}{\mathbb{R}}
\newcommand{\Z}{\mathbb{Z}}
\newcommand{\C}{\mathbb{C}}
\newcommand{\N}{\mathbb{N}}

\usepackage{graphicx}
\usepackage{amscd}
\usepackage{pstricks}
\usepackage{pst-plot}
\usepackage{multido}
\usepackage{pst-coil}
\usepackage{colortab}

\usepackage{psfrag} 

\definecolor{azul1}{rgb}{0.00,0.25,0.50}
\definecolor{azul2}{rgb}{0.00,0.50,0.75}
\definecolor{azul3}{rgb}{0.50,0.50,1.00}
\definecolor{azul4}{rgb}{0.66,0.80,0.95}

\definecolor{darkgreen}{rgb}{0,.5,0}
\definecolor{verde}{rgb}{0.55,0.68,0.09}
\definecolor{pistacho}{rgb}{0.80,0.89,0.17}
\definecolor{verdeclaro}{rgb}{0.9,0.9,0.25}
\definecolor{marron}{rgb}{0.66,0.47,0.13}
\definecolor{rojo}{rgb}{0.68,0.13,0.04}
\definecolor{azul}{rgb}{0.38,0.53,0.72}
\definecolor{amarillo}{rgb}{0.85,0.81,0.00}
\definecolor{marronclaro}{rgb}{.749,.675,.376}
\definecolor{palegray}{rgb}{0.85,0.85,0.85}
\definecolor{gray}{rgb}{0.63,0.63,0.63}
\definecolor{darkgray}{rgb}{0.50,0.50,0.50}

\theoremstyle{plain}
\newtheorem{theorem}{\bf Theorem}[section]

\theoremstyle{remark}
\newtheorem{remark}[theorem]{\bf Remark}
\newtheorem{definition}[theorem]{\bf Definition}

\begin{document}
\title{Live load matrix recovery from scattering data in linear elasticity }
\author{Juan Antonio Barcel\'{o}, Carlos Castro, 
 Mari Cruz Vilela}
\thanks{Supported by the Spanish Grant MTM2017-85934-C3-3-P.}
\begin{abstract}
We study the numerical approximation of the inverse scattering problem in the two-dimensional homogeneous isotropic linear elasticity with an unknown linear load given by a square matrix.
For both backscattering data and fixed-angle scattering data, we show how to obtain numerical approximations of the so-called Born approximations and propose new iterative algorithms that provide sequences of approximations to the unknown load. Numerical evidences of the convergence for not too large loads are also given. 
\end{abstract}
\maketitle
\markright{Live load matrix recovery from scattering data in linear elasticity}
\markleft{J.A. Barcel\'o, C. Castro and M.C. Vilela.}
\section{Introduction and statement of results}

The propagation of time-harmonic elastic waves in a two-dimensional homogeneous and isotropic
linear elastic medium is governed by the equation
\begin{equation}
\Delta^{\ast}\mathbf{u}(x)+\omega^{2}\mathbf{u}(x) = \mathrm{Q}(x)\mathbf{u}%
(x), \qquad\qquad \omega>0,\ x\in\mathbb{R}^{2}, 
\label{ecuacionQ}%
\end{equation}
where $\omega$ is the frequency of the wave, $\mathbf{u},$ the displacement vector, is a vector-valued function from
$\mathbb{R}^{2}$ to $\mathbb{R}^{2}$, and
$\Delta^*$ is  the \emph{Lamé operator} given by
\begin{equation}
\label{operador}\Delta^{\ast}\mathbf{u}(x) = \mu\Delta\mathrm{I}%
\mathbf{u}(x)+(\lambda+\mu)\nabla div\,\mathbf{u}(x),
\end{equation}
with $\Delta\mathrm{I}$ denoting the diagonal matrix with the Laplace operator
on the diagonal. 
Here $\lambda$ and $\mu$ are constants, known as the Lam\'{e} constants, that we assume to satisfy $\mu>0$ and $2\mu+\lambda>0$, so that $\Delta^*$ is  a strongly elliptic operator, and $\mathrm{Q}$ is a square real matrix of order $2$ compactly supported with support in $B(0, R)$ where $R>0$. 
The matrix $\mathrm{Q}$ describes either a linear distributed load on the elastic medium or the mass density. In the latter case $\mathrm{Q}(x)=\omega^2 (1-m(x))\mathrm{I}$ with $m(x)$ being the density and $\mathrm{I}$ denoting the identity matrix of order 2.

We are interested in recovering numerically the unknown matrix load $\mathrm{Q}$ from scattering data associated to some incident planar waves. In general, such inverse scattering problems for elastic waves have significant applications in diverse scientific areas such as geophysical exploration, nondestructive testing  or medical diagnostics \cite{Landau,UhlmannWW}.   

From the theoretical point of view, it is still an open problem whether the scattering data allow to recover the load $\mathrm{Q}$. One expects to recover $\mathrm{Q}$ from high frequency limits of
scattering data, since this is the case for the Schr\"odinger equation, the scalar analogous of \eqref{ecuacionQ}, see \cite{PSe}. As far as we know, in the case of the Lam\'e system, the recovery from high frequency limits works for
\eqref{ecuacionQ} only in very special cases, namely either assuming
$\mathrm{Q}=q(x)\mathrm{I}$ or $\mu+\lambda=0$ (see \cite[Corolary 3.5 and 3.6]{bfprv1}). From the practical point of view, such high 
frequency limits introduce an important numerical difficulty due to the possible interaction between the high frequency and the necessary discretization parameter.  

Alternatively, some information can be obtained from suitable linearizations  that do not require such high frequency limit. These are called \emph{Born approximations}, and they were defined in \cite{bfprv3} and \cite{bfprv1} according to different types of scattering data, backscattering or fixed angle respectively, with the aim of recovering singularities in the matrix load. The sharper results in recovering singularities for the scalar case (Helmholtz equation) can be found in \cite{M1} and \cite{M3}. This may suggest that the results in \cite{bfprv3} and \cite{bfprv1} can be improved. 

In this work we approximate numerically these Born approximations and show that they are in fact a very good approximation of the matrix load as long as it is not too large. We also propose iterative algorithms to recover this load and show some numerical evidence of their convergence, again when $\mathrm{Q}$ is not large.   

Let us describe more precisely the problem. First of all, to define the scattering data we introduce the usual asymptotic boundary conditions that guarantee the uniqueness of solutions for \eqref{ecuacionQ}. It is well known that when $\mathrm{Q}=0$ planar wave solutions can be divided into \textit{longitudinal waves or pressure waves} (\textit{p-waves}) and \textit{transverse waves or shear waves} (\textit{s-waves}) that propagate with different speeds in the elastic medium.
We will denote by $k_{p}$ and $k_{s}$ the speeds of propagation of p-waves and s-waves respectively, which
are given by
\begin{equation*}
\displaystyle k_{p}^{2}=\frac{\omega^{2}}{(2\mu+\lambda)}
\qquad\mathtt{\mathrm{and}} \qquad k_{s}^{2}=\frac{\omega^{2}}{\mu}.
\end{equation*}

Now, given an incident planar wave $\mathbf{u}_i$, if $\mathrm{Q}$ is a symmetric matrix with each component belonging to the space $C^1(\mathbb{R}^2)$ and compactly supported with support in $B(0, R)$ where $R>0$ then,
there exists an unique solution 
$\mathbf{u}\in \mathbf{H}^1_{loc}(\mathbb{R}^2)
=H^1_{loc}(\mathbb{R}^2)\times H^1_{loc}(\mathbb{R}^2)$ of the equation (\ref{ecuacionQ}), called \emph{the scattering solution},
which can be written as 
\begin{equation}
\mathbf{u}=\mathbf{u_{i}}+\mathbf{v}, \label{uuiv}%
\end{equation}
with $\mathbf{v},$ the \emph{scattered wave} satisfying the \emph{outgoing Kupradze radiation conditions}. They are given by
\begin{eqnarray}
&
(\partial_{r}-ik_{p})\mathbf{v}_{p}=\mathbf{o}(r^{-1/2}),\qquad
r=|x|\rightarrow\infty,\label{radiacionup}\\
&
(\partial_{r}-ik_{s})\mathbf{v}_{s}=\mathbf{o}(r^{-1/2}),\qquad
r=|x|\rightarrow\infty, \label{radiacionus}%
\end{eqnarray}
where $\mathbf{v}_{p}$ and $\mathbf{v}_{s}$ are the compressional and the shear part of $\mathbf{v}$ and are defined outside of $B(0,R)$ by
\begin{equation}
\label{vpvs}
\displaystyle\mathbf{v}_{p}=-\frac{1}{k_{p}^{2}}\nabla \textrm{ div }\mathbf{v}
\qquad
\mathtt{\mathrm{and}}
\qquad\mathbf{v}_{s}=\mathbf{v}-\mathbf{v}_{p}.
\end{equation}
See \cite{kupradze} and \cite{ikehata} for details, and also \cite[Theorem 4.1]{bfprv5} for the three-dimensional case.
Existence and uniqueness results with less regularity in the components of $\mathrm{Q}$ are only partial and for sufficiently large $\omega$  (see \cite[Proposition 3.1]{bfprv1} for details).  

Scattering data are now obtained from the asymptotic behavior as $|x|\to \infty$ of scattered waves with particular incident waves. Observe that these scattered waves $\mathbf{v}=\mathbf{v}(\cdot;\omega,\mathbf{u}_{i})$ are perturbations of some incident waves due to the load $\mathrm{Q}$, see (\ref{uuiv}), and  they satisfy the following equation:
\begin{equation}
\label{ecuacionv}
(\Delta^{\ast}+\omega^{2}\mathrm{I})\mathbf{v}
=
\mathrm{Q}(\mathbf{u}_i+\mathbf{v})
\end{equation}
together with the outgoing Kupradze conditions \eqref{radiacionup} and \eqref{radiacionus}. Let us denote by $\mathbf{f}$ the second hand term in \eqref{ecuacionv}. It is well known that, when $\mathbf{f}\in\mathbf{C}_0^\infty(\mathbb{R}^2)$, the outgoing solution of the equation  \eqref{ecuacionv} can be written outside the support of $\mathbf{f}$ as $\mathbf{v}=\mathbf{v}_p+\mathbf{v}_s$, with $\mathbf{v}_p$ and $\mathbf{v}_s$ given in \eqref{vpvs}. 
One can check that $\mathbf{v}_p$ and $\mathbf{v}_s$ are solutions of the following vectorial Helmholtz equations (see the proof of Theorem 1.1 in \cite[pp.833-834]{bfprv2})
\begin{align*}
\Delta\mathrm{I}\mathbf{v}_{p}+k_{p}^{2}\mathbf{v}_{p}  &  =-\frac{1}%
{\omega^{2}}\nabla\operatorname{div}\mathbf{f},\\
\Delta\mathrm{I}\mathbf{v}_{s}+k_{s}^{2}\mathbf{v}_{s}  &  =\frac{1}%
{\omega^{2}}(k_{s}^{2}\mathbf{f}+\nabla\operatorname{div}\mathbf{f}).
\end{align*}
Since $\mathbf{v}_p$ and $\mathbf{v}_s$ also satisfy the corresponding outgoing Sommerfeld radiation conditions (see \eqref{radiacionup} and \eqref{radiacionus}) it is easy to obtain the following asymptotic expansions
\begin{align} \label{eq:asy1}
\mathbf{v}_{p}(x)  &  = c\ k_{p}^{-1/2}\,\frac{e^{ik_{p}|x|}}%
{|x|^{1/2}}\,\mathbf{v}_{p,\infty}\left(  k_{p},x/|x|\right)  +
\mathbf{o}\left(  |x|^{-1/2}\right)  , \qquad|x|\rightarrow\infty,\\
\mathbf{v}_{s}(x)  &  = c\ k_{s}^{-1/2}\,\frac{e^{ik_{s}|x|}}%
{|x|^{1/2}}\,\mathbf{v}_{s,\infty}\left(  k_{s},x/|x|\right)  +
\mathbf{o}\left(  |x|^{-1/2}\right)  , \qquad|x|\rightarrow\infty.
\label{eq:asy2}
\end{align} 
for some $\mathbf{v}_{p,\infty}$ and $\mathbf{v}_{s,\infty}$ that are known as the longitudinal and transverse scattering amplitudes of $\mathbf{v}$ respectively.
Moreover, they can be written as
\begin{align*}
\mathbf{v}_{p,\infty}\left(  k_{p},x/|x|\right)   
&  
=-\frac{1}{\omega^{2}%
}\,(\nabla div\,\mathbf{f})\ \widehat{}\,\left(  k_{p}x/|x|\right),
\\
\mathbf{v}_{s,\infty}\left(  k_{s},x/|x|\right)   
&  
=\frac{1}{\omega^{2}%
}((\nabla div+k_{s}^{2}\mathrm{I})\mathbf{f})\ \widehat{}\,\left(
k_{s}x/|x|\right). 
\end{align*}
These definitions were extended in \cite{bfprv1} to the case $\mathbf{f}\in\mathbf{L}^p(\mathbb{R}^2)$ for some $1\leq p<6/5$ by using the 
Leray projection $\mathrm{I}-\mathcal{R}$, defined via the Fourier transform  by 
\begin{equation}
\label{leray}
(\mathcal{R}\mathbf{f})\,\widehat{\ }\,(\xi)
=\left(\widehat{\mathbf{f}}(\xi)\cdot\frac{\xi}{|\xi|}\right)\frac{\xi}{|\xi|},\qquad \xi\in\mathbb{R}^2.
\end{equation}
More precisely, in such a case, the scattering amplitudes of the solution $\mathbf{v}$ are defined by  (see \cite[Definition 2.3]{bfprv1})
\begin{align*}
\mathbf{v}_{p,\infty}\left(  k_{p},x/|x|\right)  
&=\frac{1}{2\mu+\lambda
}\,(\mathcal{R}\mathbf{f})\ \widehat{}\,\left(  k_{p}x/|x|\right)  ,
\\
\mathbf{v}_{s,\infty}\left(  k_{s},x/|x|\right)  
&=\frac{1}{\mu}((\mathrm{I}%
-\mathcal{R})\mathbf{f})\ \widehat{}\,\left(  k_{s}x/|x|\right)  .
\end{align*}
Note that, under the conditions that ensure the existence and uniqueness of the direct scattering problem, we have that the right hand side of  \eqref{ecuacionv} belongs to $\mathbf{L}^p(\mathbb{R}^2)$ with $1\leq p<6/5$. Therefore the scattering amplitudes of the scattered solution $\mathbf{v}$ are well defined.

Moreover, writing $(\mathcal{R}\mathbf{f})\,\widehat{\ }\,(\xi)=\Pi_{\xi/|\xi| }\widehat{ \mathbf{f}}(\xi)$, we have that the \emph{scattering amplitudes} of $\mathbf{v}$ are given by
\begin{align}
\mathbf{v}_{p,\infty}\left(x/|x|;\omega,\mathbf{u}_i\right)  
&=\frac{1}{2\mu+\lambda}\,
\Pi_{x/|x| }(\mathrm{Q}(\mathbf{u}_i+\mathbf{v}))\, \widehat{}\,\left(  k_{p}x/|x|\right)  ,
\label{amplitudvp}
\\
\mathbf{v}_{s,\infty}\left(x/|x|;\omega,\mathbf{u}_i\right)  
&=\frac{1}{\mu}\,
(\mathrm{I}-\Pi_{x/|x| })
(\mathrm{Q}(\mathbf{u}_i+\mathbf{v}))\, \widehat{}\,\left(  k_{s}x/|x|\right)  .
\label{amplitudvs}
\end{align}
              
On the other hand, as incident waves we will always consider plane waves, either 
longitudinal plane waves (\textit{plane p-waves})
\[
\label{uip}\mathbf{u}^{p}_{i}(x)=e^{ik_{p}\theta\cdot x}\theta,
\]
where the wave direction and the polarization vector are the same $\theta\in\mathbb{S}^{1}$
or,
transverse plane waves
(\textit{plane s-waves}) for which the the polarization vector is orthogonal to the wave direction, that is
\[
\label{uis}\mathbf{u}^{s}_{i}(x)=e^{ik_{s}\theta\cdot x}\theta^\perp,
\]
where $\theta^\perp$ is a unitary vector orthogonal to $\theta$.
For the corresponding scattered waves, denoted by $\mathbf{v}^p$ and $\mathbf{v}^s$ respectively, we define different scattering amplitudes, which will depend on
the corresponding parameters, namely 
$\mathbf{v}_{p,\infty}^{p}(\zeta;\omega,\theta),$ 
$\mathbf{v}_{s,\infty}^{p}(\zeta;\omega,\theta),$
$\mathbf{v}_{p,\infty}^{s}(\zeta;\omega,\theta,\theta^\perp)$ and, 
$\mathbf{v}_{s,\infty}^{s}(\zeta;\omega,\theta,\theta^\perp),$ 
where $\zeta= \frac x{|x|}.$
These data are known as $p\to p,$ $p\to s,$ $s\to p$ and, $s\to s$ scattering
data, respectively (see Section \ref{SecAmplitudes} for explicit definitions).

We have written $\mathbf{v}^s_{p,\infty}$ and $\mathbf{v}^s_{s,\infty}$ depending on the parameters $\theta$ and $\theta^\perp$. However, in dimension 2 both vectors depends on a single parameter $\theta$, except for a minus sign, since there is only one orthonormal vector to a given $\theta$. Thus, we can say that these data depend only on $\zeta$, $\omega$ and $\theta$, up to a sign. 

In the \emph{inverse elastic scattering problem} that we study, $\mathrm{Q}$ is assumed to be an unknown matrix that we determine from the knowledge of these scattering amplitudes. The Born approximation is formally obtained by considering $\textbf{v}=0$ in \eqref{amplitudvp}-\eqref{amplitudvs}. Note  that there is no guarantee that, under this hypotheses, system  \eqref{amplitudvp}-\eqref{amplitudvs} has a solution. In fact, when considering the particular class of aforementioned incident planar waves this system becomes overdetermined since we have 4 scattering amplitudes, which are vector-valued functions depending on 3 parameters, to determine $\mathrm{Q}$ that is characterized by 4 scalar functions depending on 2 variables. It is then natural to restrict the scattering data. The most commonly used partial scattering data are the \emph{backscattering data} and the \emph{fixed-angle scattering data}. In the first case, the receptor direction is opposite to the incident direction, that is, $\zeta=-\theta$, and in the latter, the incident direction $\theta$ is fixed. 

The strategy to clear the four components in the Born approximation of $\mathrm{Q}$ from \eqref{amplitudvp}-\eqref{amplitudvs} with $\textbf{v}=0$ is not easy as described in Section 3 below. Moreover, it is remarkable that the definition of the Born approximation for backscattering data
involves a linear combination of the $p\to p,$ $p\to s,$ $s\to p$ and $s\to s$ scattering data (see Definition \ref{BornB} and Remark \ref{datosB} bellow), while in the definition for the fixed-angle data only the $p\to p$ and $p\to s$ or the $s\to p$ and $s\to s$ data are needed, depending on whether $k_s\ge k_p$ (see Definitions \ref{BornAK>1} and \ref{BornAK<1} bellow). This means that in the backscattering we need both, longitudinal and transverse plane waves, while in the fixed-angle data we only need one of them, the one with lower frequency (see Remark \ref{data} bellow).

We show a trigonometric collocation method to approximate numerically the Born approximations for both backscattering data and fixed-angle scattering data. We also propose, in each case, iterative algorithms that improve the Born approximation to recover the unknown load $\mathrm{Q}$ from its scattering data. Roughly speaking, if we denote by $\mathrm{Q}_B$ the Born approximation obtained from a class of scattering data, these algorithms are based on the following iterative process:

\medskip

{\bf Algorithm 1}
\begin{itemize}
    \item[$\star$] 
    Let $\mathrm{Q}_1 =\chi_R \mathrm{Q}_B$.
    \\[0.5ex]
    \item[$\star$]
    For $n\ge 1,$ compute $\mathrm{Q}_{n+1}$ from $\mathrm{Q}_{n}$ as follows,
    \begin{enumerate} 
    \item  solve \eqref{ecuacionv} with
    $\mathrm{Q}=\mathrm{Q}_{n}$
    to get 
    $\mathbf{v}_{n},$
    \item construct $\widetilde{\mathrm{Q}}_{n+1}$ from $\mathbf{v}_{n}$
    using,
    \\[1ex]
    $
    \displaystyle
    \mathbf{v}_{p,\infty}\left(x/|x|;\omega,\mathbf{u}_i\right)  
    =\frac{1}{2\mu+\lambda}\,
    \Pi_{x/|x| }( \widetilde{\mathrm{Q}}_{n+1}\mathbf{u}_i+\mathrm{Q}_{n}\mathbf{v}_n))\, \widehat{}\,\left(  k_{p}x/|x|\right)  ,
    $
    \\[1ex]
    $
    \displaystyle
    \mathbf{v}_{s,\infty}\left(x/|x|;\omega,\mathbf{u}_i\right)  
    =\frac{1}{\mu}\,
    (\mathrm{I}-\Pi_{x/|x| })
    (\widetilde{\mathrm{Q}}_{n+1}\mathbf{u}_i+\mathrm{Q}_{n}\mathbf{v}_n)\, \widehat{}\,\left(  k_{s}x/|x|\right),
    $
    \item set $\mathrm{Q}_{n+1}=\chi_R\widetilde{\mathrm{Q}}_{n+1}$.
    \end{enumerate}
\end{itemize}
Here, $\chi_R$ is the characteristic function of the ball  $B(0,R)$ that is necessary since the Born approximations may not be of compact support, while $\mathrm{Q}$ does.
As in the Born approximation, the matrix $\widetilde{\mathrm{Q}}_{n+1}$ is obtained combining particular, and complicated, choices of planar incident incident waves $\mathbf{u}_i$ that depend on the scattering data. In Section 5 below we give a more detailed version of this algorithm.

We would like to point out that $\mathbf{v}_n$ denotes the scattered wave associated to the load $\mathrm{Q}_n$, and to get it in each iteration it is necessary to solve \eqref{ecuacionv} with $\mathrm{Q}=\mathrm{Q}_{n}$.
As a consequence, in each iteration an integral equation has to be solved, and this makes the algorithm computationally cost. The idea behind this algorithm was first introduced in \cite{BCR} to recover a potential in the inverse quantum scattering problem, which can be seen as a scalar version of the elasticity problem in certain very simple cases,
for similar methods see \cite{N} and \cite{AHN}.
Later, in \cite{BCLV}, for the inverse quantum scattering problem, a new algorithm was constructed that avoids the solution of integral equations, but it is not clear how to extend this to the elasticity system, which is not a scalar problem.

Numerical experiments show that the sequences provided by the constructed algorithms converge rapidly to the matrix $\mathrm{Q}$
(for more details see Section \ref{experimentos}).
The study of this convergence from an analytical point of view is very interesting, but also hard, so we leave it for future work. Note that the convergence of the algorithms to the matrix load $\mathrm{Q}$ would provide, in particular, uniqueness of the inverse scattering problem with the associated scattering data, i.e backscattering or fixed angle.

Both the Born approximations and the algorithms proposed here are also valid in higher dimensions (see Remarks \ref{dimensionB} and \ref{dimensionAF}). However, the computational cost in dimension 3 increases significantly.

The rest of the article is organized as follows. 
In Section 2 we explicitly write down the scattering amplitudes from which we will obtain the scattering data we will use to solve the two inverse problems we are interested in, namely the backscattering and the fixed-angle scattering problem.
In Section 3 we use these scattering data to define the corresponding Born approximations.
In Section 4 we explain the algorithms that will allow us to numerically recover the unknown matrix $\mathrm{Q}$ from its Born approximations.
In Section \ref{Numerico} we describe the numerical method to construct the Born approximations and the sequences of approximations to the load $\mathrm{Q}$.
Finally, in Section \ref{experimentos} we present numerical experiments.
\section{The scattering amplitudes.}
\label{SecAmplitudes}
In this section we give precise definitions of the scattering amplitudes that we use later to obtain the scattering data according to the different strategies to obtain the Born approximation, namely backscattering or fixed-angle scattering data.  

For convenience, we introduce the outgoing resolvent of the operator $\Delta^*$ provided by the limiting absorption principle given in Theorem 1.1 of \cite{bfprv2} and denoted by $\mathbf{R}(\omega^2+i0)$.
If we apply this operator to \eqref{ecuacionv}, we have that $\mathbf{v}$ satisfies the following Lippmann-Schwinger equation:
\begin{equation}
\label{equacionintegral}
\mathbf{v}
=
\mathbf{R}(\omega^2+i0)\mathrm{Q}(\mathbf{u}_i+\mathbf{v}).
\end{equation}

Therefore, if we consider an incident plain p-wave $\mathbf{u}_{i}%
^{p}(x)=e^{ik_{p}\theta\cdot x}\theta$ with $\theta\in\mathbb{S}^1,$ 
the scattered wave, $\mathbf{v}^{p}(\cdot;\omega,\theta)$,
satisfies the following integral equation:
\begin{equation}
\label{equacionintegralvp}
\mathbf{v}^p(x;\omega,\theta)
=
\mathbf{R}(\omega^2+i0)\left(\mathrm{Q}e^{ik_{p}\theta\cdot (\cdot)}\theta \right)(x)
+
\mathbf{R}(\omega^2+i0)\left(\mathrm{Q}\mathbf{v}^{p}(\cdot;\omega,\theta) \right)(x).
\end{equation}
Moreover, from \eqref{amplitudvp} and \eqref{amplitudvs}
we obtain that the $p\rightarrow p$ scattering data is given by
\begin{equation}
\mathbf{v}_{p,\infty}^{p}(\zeta;\omega,\theta)
=\frac{1}{2\mu+\lambda}\,
\left(
\Pi_{\zeta} \widehat{\mathrm{Q}}\left(  k_{p}(\zeta-\theta)\right)
\theta+
\Pi_{\zeta} (\mathrm{Q}\mathbf{v}^{p}%
(\cdot;\omega,\theta))\ \widehat{}\ \left(  k_{p}\,\zeta\right) 
\right) ,
\label{amplitudpp}%
\end{equation}
and the $p\rightarrow s$ scattering data by
\begin{equation}
\mathbf{v}_{s,\infty}^{p}(\zeta;\omega,\theta)
=\frac{1}{\mu}\,
\left(
(\mathrm{I}-\Pi_{\zeta}) \widehat{\mathrm{Q}}\left(  k_{s}\zeta-k_{p}\theta\right)
\theta+ 
(\mathrm{I}-\Pi_{\zeta})(\mathrm{Q}\mathbf{v}^{p}%
(\cdot;\omega,\theta))\ \widehat{}\ \left(  k_{s}\zeta\right) 
\right) ,
\label{amplitudps}%
\end{equation}
where $\zeta=x/|x|$.

Whereas if we consider as the incident wave a plane s-wave, 
$\mathbf{u}_i^s(x)=e^{ik_{s}\varphi\cdot x}\phi$ where $\varphi$ and $\phi$ are unitary
orthogonal vectors, the scattered wave $\mathbf{v}^{s}(\cdot;\omega,\varphi,\phi)$ satisfies the following integral equation: 
\begin{equation}
\label{equacionintegralvs}
\mathbf{v}^s(x;\omega,\varphi,\phi)
=
\mathbf{R}(\omega^2+i0)\left(\mathrm{Q}e^{ik_{s}\varphi\cdot (\cdot)}\phi \right)(x)
+
\mathbf{R}(\omega^2+i0)\left(\mathrm{Q}\mathbf{v}^{s}(\cdot;\omega,\varphi,\phi) \right)(x).
\end{equation}
And from \eqref{amplitudvp} and \eqref{amplitudvs} we have that the
$s\rightarrow p$ scattering data is given by
\begin{equation}
\label{amplitudsp} 
\mathbf{v}_{p,\infty}^{s}(\zeta;\omega,\varphi,\phi)    
=\frac{1}{2\mu+\lambda}\,
\left(
\Pi_{\zeta}\widehat{\mathrm{Q}}(k_{p}\zeta-k_{s}\varphi)\phi 
+
\Pi_{\zeta}(\mathrm{Q}\mathbf{v}^{s}(\cdot
;\omega,\varphi,\phi))\ \widehat{}\ (k_{p}\zeta)
\right),
\end{equation}
and  the $s\rightarrow s$ scattering data by
\begin{equation}
 \label{amplitudss} 
\mathbf{v}_{s,\infty}^{s}(\zeta;\omega,\varphi,\phi)    
=\frac{1}{\mu}
\left(
(\mathrm{I}-\Pi_{\zeta}) 
\widehat{\mathrm{Q}}\left(k_{s}(\zeta-\varphi)\right)\phi
+
(\mathrm{I}-\Pi_{\zeta}) (\mathrm{Q}\mathbf{v}^{s}%
(\cdot;\omega,\varphi,\phi))\ \widehat{}\, \left(  k_{s}\zeta\right)
\right).
\end{equation}
For convenience, we introduce the parameter
\begin{equation}
\label{K}
K=\frac{k_s}{k_p}=\frac{\sqrt{2\mu+\lambda}}{\sqrt{\mu}},
\end{equation}
and we rescale (\ref{amplitudpp}), (\ref{amplitudps}), (\ref{amplitudsp}) and (\ref{amplitudss}) writing
\begin{align}
\label{amplitudppR}
(2\mu+\lambda)\mathbf{v}_{p,\infty}^{p}(\zeta;\sqrt{2\mu+\lambda}\omega,\theta)
=
&
\Pi_{\zeta}\widehat{\mathrm{Q}}\left(\omega(\zeta-\theta)\right)\theta
\\
&
+
\Pi_{\zeta}(\mathrm{Q}\mathbf{v}^{p}(\cdot;\sqrt{2\mu+\lambda}\omega,\theta))\ \widehat{}\ \left(\omega\,\zeta\right),
\nonumber
\end{align}
\begin{align}
\label{amplitudpsR}
\mu\mathbf{v}_{s,\infty}^{p}(\zeta;\sqrt{2\mu+\lambda}\omega,\theta)
=
&
(\mathrm{I}-\Pi_{\zeta})\widehat{\mathrm{Q}}\left(\omega(K\zeta-\theta)\right)\theta
\\
&
+
(\mathrm{I}-\Pi_{\zeta})(\mathrm{Q}\mathbf{v}^{p}(\cdot;\sqrt{2\mu+\lambda}\omega,\theta))\ \widehat{}\ \left(K\omega\zeta\right),
\nonumber
\end{align}
\begin{align}
\label{amplitudspR}
(2\mu+\lambda)\mathbf{v}_{p,\infty}^{s}(\zeta;\sqrt{\mu}\omega,\varphi,\phi)
=
&
\Pi_{\zeta}\widehat{\mathrm{Q}}\left(\omega(K^{-1}\zeta-\varphi)\right)\phi
\\
&
+
\Pi_{\zeta}(\mathrm{Q}\mathbf{v}^{s}(\cdot;\sqrt{\mu}\omega,\varphi,\phi))\ \widehat{}\ \left(K^{-1}\omega\,\zeta\right),
\nonumber
\end{align}
\begin{align}
\label{amplitudssR}
\mu\mathbf{v}_{s,\infty}^{s}(\zeta;\sqrt{\mu}\omega,\varphi,\phi)
=
&
(\mathrm{I}-\Pi_{\zeta})\widehat{\mathrm{Q}}\left(\omega(\zeta-\varphi)\right)\phi
\\
&
+
(\mathrm{I}-\Pi_{\zeta})(\mathrm{Q}\mathbf{v}^{s}(\cdot;\sqrt{\mu}\omega,\varphi,\phi))\ \widehat{}\ \left(\omega\zeta\right).
\nonumber
\end{align}
We would like to note that $\mathbf{v}^{p}(x;\sqrt{2\mu+\lambda}\omega,\theta)$ is the scattered solution of the direct scattering problem with energy $(2\mu+\lambda)\omega^2$ and incident wave $e^{i\omega\theta\cdot x}\theta$, and $\mathbf{v}^{s}(x;\sqrt{\mu}\omega,\varphi,\phi)$ is the scattered solution of the problem with energy $\mu\omega^2$ and incident wave $e^{i\omega\varphi\cdot x}\phi$.
Moreover, they are solutions of the Lippmann-Schwinger equation
\begin{equation}
\label{LSeq}
\mathbf{v}(x)
=
\mathbf{R}(c^2+i0)\left(\mathrm{Q}\mathbf{u}_i \right)(x)
+
\mathbf{R}(c^2+i0)\left(\mathrm{Q}\mathbf{v}(\cdot) \right)(x).
\end{equation}
with
$\mathbf{u}_i(x)=e^{i\omega\theta\cdot x}\theta$ and $c^2=(2\mu+\lambda)\omega^2$, and
$\mathbf{u}_i(x)=e^{i\omega\varphi\cdot x}\phi$ and $c^2=\mu\omega^2$ respectively. 
\section{Born approximation for scattering  data.}
In this section we will use the scattering amplitudes given in \eqref{amplitudppR}--\eqref{amplitudssR} to construct good approximations for the matrix $\mathrm{Q}(x)$ called Born approximations. As we mentioned in the introduction the strategy is different for backscattering and fixed-angel scattering data and therefore we divide this section in two more subsections where we analyze both cases separately.

\subsection{Born approximation for backscattering data.}
Following \cite[Subsection 2.1]{bfprv3}, we look for a matrix, called the Born approximation for backscattering data, which is a \lq\lq good approximation\rq\rq for the matrix $\mathrm{Q}(x)$ and it is written in terms of the scattering amplitudes. Since these are vector valued functions written in terms of the Fourier transform of $\mathrm{Q}$ acting on a vector, we introduce $\{e_1,e_2\}$, the canonical base of $\mathbb{R}^2$, and for $x\in\mathbb{R}^2$ and $i=1,2,$ we write the vector $\mathrm{Q}(x)e_i$ via its Fourier transform, $\widehat{\mathrm{Q}}(\xi)e_i,$ and we make the change of variable $\xi=-2\omega\theta,$ so that we get
\begin{equation*}
	\label{QB}
	\mathrm{Q}(x)e_i 
	=\int_{\mathbb{R}^2}e^{i \xi \cdot x}\widehat{\mathrm{Q}}(\xi)e_i\,d\xi
	=2 \int_0^\infty \omega \int_{\mathbb{S}^1}e^{-2i\omega \theta \cdot x}  
	\widehat{\mathrm{Q}}(-2 \omega \theta)e_i\,  
	d \sigma(\theta)\, d \omega.
\end{equation*}
For any $\theta \in \mathbb{S}^1$ we can write
$$e_i=(e_i \cdot \theta) \theta +(e_i \cdot \theta^\perp)\theta^\perp ,\qquad i=1,2,$$
where $\theta^\perp$ is a unitary vector orthogonal to $\theta$, and thus we have that
$$
\widehat{\mathrm{Q}}(-2 \omega \theta)e_i=
(e_i \cdot \theta)\widehat{\mathrm{Q}}(-2\omega \theta) \theta 
+  ( e_i \cdot \theta^{\perp})   \widehat{\mathrm{Q}}(-2 \omega \theta) \theta^\perp.     
$$
Therefore, it is enough to write the vectors 
$\widehat{\mathrm{Q}}(-2\omega \theta) \theta$
and
$\widehat{\mathrm{Q}}(-2 \omega \theta) \theta^\perp$
in terms of the scattering amplitudes.

If we consider as incident p-wave 
$\mathbf{u}_{i}^{p}(x)=e^{ik_{p}\theta\cdot x}\theta$
and as incident s-wave
$\mathbf{u}_i^s(x)=e^{ik_{s}\theta\cdot x}\theta^\perp$, from 
\eqref{amplitudppR}--\eqref{amplitudssR} with $\zeta=-\theta$, $\varphi=\theta$ and $\phi=\theta^\perp$, we get
\begin{align}
\label{amplitudppB}
\Pi_{-\theta}\widehat{\mathrm{Q}}\left(-2\omega\theta\right)\theta
=
&
(2\mu+\lambda)\mathbf{v}_{p,\infty}^{p}(-\theta;\sqrt{2\mu+\lambda}\omega,\theta)
\\
&
-
\Pi_{-\theta}(\mathrm{Q}\mathbf{v}^{p}(\cdot;\sqrt{2\mu+\lambda}\omega,\theta))\ \widehat{}\ \left(-\omega\,\theta\right),
\nonumber
\end{align}
\begin{align}
\label{amplitudpsB}
(\mathrm{I}-\Pi_{-\theta})\widehat{\mathrm{Q}}\left(-(K+1)\omega\theta\right)\theta
=
&
\mu\mathbf{v}_{s,\infty}^{p}(-\theta;\sqrt{2\mu+\lambda}\omega,\theta)
\\
&
-
(\mathrm{I}-\Pi_{-\theta})(\mathrm{Q}\mathbf{v}^{p}(\cdot;\sqrt{2\mu+\lambda}\omega,\theta))\ \widehat{}\ \left(-K\omega\theta\right),
\nonumber
\end{align}
\begin{align}
\label{amplitudspB}
\Pi_{-\theta}\widehat{\mathrm{Q}}\left(-(K^{-1}+1)\omega\theta\right)\theta^\perp
=
&
(2\mu+\lambda)\mathbf{v}_{p,\infty}^{s}(-\theta;\sqrt{\mu}\omega,\theta,\theta^\perp)
\\
&
-
\Pi_{-\theta}(\mathrm{Q}\mathbf{v}^{s}(\cdot;\sqrt{\mu}\omega,\theta,\theta^\perp))\ \widehat{}\ \left(-K^{-1}\omega\,\theta\right),
\nonumber
\end{align}
\begin{align}
\label{amplitudssB}
(\mathrm{I}-\Pi_{-\theta})\widehat{\mathrm{Q}}\left(-2\omega\theta\right)\theta^\perp
=
&
\mu\mathbf{v}_{s,\infty}^{s}(-\theta;\sqrt{\mu}\omega,\theta,\theta^\perp)
\\
&
-
(\mathrm{I}-\Pi_{-\theta})(\mathrm{Q}\mathbf{v}^{s}(\cdot;\sqrt{\mu}\omega,\theta,\theta^\perp))\ \widehat{}\ \left(-\omega\theta\right).
\nonumber
\end{align}

From \eqref{amplitudppB} and a rescaled version of \eqref{amplitudpsB} where $\omega$ is replaced by $2\omega/(K+1)$ we obtain that
\begin{equation}
\label{parte1B}
\widehat{\mathrm{Q}}\left(-2\omega\theta\right)\theta=
\mathbf{v}^p_\infty(\omega,\theta)-\mathbf{e}^p(\omega,\theta),
\end{equation}
where $\mathbf{v}^p_\infty$ is written in terms of the $p\rightarrow p$ and $p\rightarrow s$  scattering data as follows
\begin{equation}
\label{b_p}
\mathbf{v}_{\infty}^{p}(\omega,\theta)
=
(2\mu+\lambda)\mathbf{v}_{p,\infty}^{p}(-\theta;\sqrt{2\mu+\lambda}\omega,\theta)
+
\mu\mathbf{v}_{s,\infty}^{p}\left(-\theta;\frac{2\sqrt{2\mu+\lambda}}{K+1}\omega,\theta\right), 
\end{equation}
and $\mathbf{e}^p$ is an error term given by
\begin{equation}%
\begin{split}
\mathbf{e}^{p}(\omega,\theta)
=  
&  
\Pi_{-\theta}\left[  (\mathrm{Q}%
\mathbf{v}^{p}(\cdot;\sqrt{2\mu+\lambda}\omega,\theta))\,\widehat{}\ \left(
-\omega\,\theta\right)  \right] 
\\
&  (\mathrm{I}-\Pi_{-\theta})\left[\left(\mathrm{Q}\mathbf{v}^{p}\left(\cdot
;\frac{2\sqrt{2\mu+\lambda}}{K+1}\omega,\theta\right)\right)\mbox{\raisebox{0.2cm}{$\widehat{}$}}\ \left(  -\frac{2K}{K+1}\omega\theta\right)
\right].
\end{split}
\label{hp}%
\end{equation}
And from \eqref{amplitudssB} and a rescaled version of \eqref{amplitudspB} where $\omega$ is replaced by $2\omega/(K^{-1}+1)$ we have that
\begin{equation}
\label{parte2B}
\widehat{\mathrm{Q}}\left(-2\omega\theta\right)\theta^\perp
=
\mathbf{v}^s_\infty(\omega,\theta)-\mathbf{e}^s(\omega,\theta),
\end{equation}
where
$\mathbf{v}^s_\infty$ is written in terms of the $s\rightarrow p$ and $s\rightarrow s$ scattering data as follows
\begin{equation}
\label{b_s}
\mathbf{v}_{\infty}^{s}(\omega,\theta) 
=   
(2\mu+\lambda)
\mathbf{v}_{p,\infty}^{s}\left(-\theta;\frac{2\sqrt{\mu}}{K^{-1}+1}\omega,\theta,\theta^\perp\right) 
+ 
\mu\mathbf{v}_{s,\infty}^{s}(-\theta;\sqrt{\mu}\omega\,\theta,\theta^\perp),
\end{equation}
and $\mathbf{e}^s$ is an error term given by
\begin{equation}%
\begin{split}
\mathbf{e}^{s}(\omega,\theta)
=  &  
\Pi_{-\theta}\left[\left(\mathrm{Q}
\mathbf{v}^{s}\left(\cdot;\frac{2\sqrt{\mu}}{K^{-1}+1}\omega,\theta,\theta^\perp\right)\right)  \mbox{\raisebox{0.2cm}{$\widehat{}$}}
\ \left(
-\frac{2K^{-1}}{K^{-1}+1}\omega\theta\right)  \right] 
\\
& 
+(\mathrm{I}-\Pi_{-\theta})\left[  (\mathrm{Q}\mathbf{v}^{s}(\cdot
;\sqrt{\mu}\omega,\theta,\theta^\perp))\ \widehat{}\ \left(  -\omega
\,\theta\right)  \right].
\end{split}
\label{hs}%
\end{equation}

Using \eqref{parte1B} and \eqref{parte2B} we can write
\begin{equation}
\label{gorroQ}
\widehat{\mathrm{Q}}(-2 \omega \theta)e_i
=
(e_i \cdot \theta)\mathbf{v}^p_\infty(\omega,\theta)
+
( e_i \cdot \theta^{\perp})\mathbf{v}^s_\infty(\omega,\theta)
-\widehat{\mathrm{E}}(-2 \omega \theta)e_i,
\end{equation}
where $\mathrm{E}=\mathrm{E}(\mathrm{Q})$ is an error term given by
\begin{equation}
\label{errorb}
\widehat{\mathrm{E}}(-2 \omega \theta)e_i
=
(e_i \cdot \theta)\mathbf{e}^p(\omega,\theta)
+
( e_i \cdot \theta^{\perp})\mathbf{e}^s(\omega,\theta),
\end{equation}
with $\mathbf{e}^p$ and $\mathbf{e}^s$ defined in \eqref{hp} and \eqref{hs} respectively.

Then, it is natural to introduce the following definition.
\begin{definition}
\label{BornB}
Given $\mathrm{Q}$ a matrix of order $2$,
we define its \it{Born approximation for backscattering data} as a matrix $\mathrm{Q}_B$ such that 
\begin{equation}
\label{Def_BornB}
{\widehat{\mathrm{Q}_{B}}}(\xi)e_{i}
=
(e_{i} \cdot \theta)\mathbf{v}_{\infty}^{p}(\omega,\theta)
+
( e_{i}\cdot \theta^\perp)\mathbf{v}_{\infty}^{s}(\omega,\theta),
\qquad
\xi\in\mathbb{R}^2,\ 
i=1,2,
\end{equation}
where $\{e_1,  e_2\}$ is the canonical base of $\mathbb{R}^{2},$ 
$\theta=\theta(\xi)\in\mathbb{S}^1$ and $\omega=\omega(\xi)>0$ are such that
$\xi=-2 \omega \theta $, and 
$\mathbf{v}_{\infty}^{p}$ and $\mathbf{v}_{\infty}^{s}$ are given in \eqref{b_p} and \eqref{b_s} respectively.
\end{definition}
From this definition and \eqref{gorroQ} we have that 
\begin{equation}
\label{definicion_juan}
\mathrm{Q}(x)=\mathrm{Q}_B(x)-\mathrm{E}(x),\qquad x\in\mathbb{R}^2,
\end{equation}
where $\mathrm{E}=\mathrm{E}(\mathrm{Q})$ is the error term defined by \eqref{errorb}.

\begin{remark}
	\label{OmegaTheta}
	Observe that any $\xi\in\mathbb{R}^2$ can be uniquely represented as $\xi=-2 \omega \theta $ with $\theta\in\mathbb{S}^1$ and $\omega>0$ given by
	\begin{equation}
	\label{omegatheta}
	\omega=\frac{|\xi|}{2}\quad \text{and}\quad \theta=-\frac{\xi}{|\xi|}.
	\end{equation}
\end{remark}
\begin{remark}
	\label{dimensionB}
	The previous definition can be extended to higher dimensions $n>2$ with a matrix $\mathrm{Q}$ (see \cite[Definition 2.7]{bfprv3}), and notice that it depends on the choice of $\theta^\perp$.
\end{remark}
\begin{remark}
	\label{datosB}
	We would like to note that in the definition of the Born approximation for backscattering data we have used four scattering data, namely 
	a $p\rightarrow p$ data associated to the incident wave $e^{i\omega\theta\cdot x}\theta$, 
	a $p\rightarrow s$ data associated to $e^{i\frac{2\omega}{K+1}\theta\cdot x}\theta$,
	a $s\rightarrow p$ data associated to $e^{i\frac{2\omega}{K^{-1}+1}\theta\cdot x}\theta^\perp$,
	and 
	a $s\rightarrow s$ data associated to $e^{i\omega\theta\cdot x}\theta^\perp$.
\end{remark}
\begin{remark}
	\label{mu+lambda=1}
	When $\mu+\lambda=0$ the Lamé operator given in \eqref{operador} becomes $\Delta\mathrm{I}$ and thus, the equation \eqref{ecuacionQ} is a vectorial Helmholtz equation. 
	Moreover, if $\mathrm{Q}=q\mathrm{I}$ with $q$ a scalar function, in \eqref{ecuacionQ} we have two uncoupled Helmhotz equations, one for each component of $\mathbf{u}$.
	More precisely, denoting this components by $u_1$ and $u_2$, we have the equations
	\begin{equation}
	\label{Helmholtzeq}
	\Delta u_j+k_s^2 u_j=\frac{q}{\mu}u_j, \qquad j=1,2.
	\end{equation}
	In this case the Born approximation for backscattering data of $\mathrm{Q}$ is $\mathrm{Q}_B=q_B\mathrm{I}$ with $q_B$ given in terms of Born approximations of $q/\mu$.
	With a few calculations one can check that writing $\xi=-2\omega\theta$ with $\omega>0$ and $\theta=(\theta_1,\theta_2)\in\mathbb{S}$, $q_B$ is given by
	\begin{equation*}
	\label{ElasticidadHelmholtz}
	\widehat{q_B}(-2\omega\theta)=\mu
	\left(\theta_1\widehat{q^1_B}(-2\omega\theta)+
	\theta_2\widehat{q^2_B}(-2\omega\theta)\right),
	\end{equation*}
	where $q^j_B$ is the Born approximation for backscattering data of $q/\mu$ with incident wave 
	$e^{ik_s\theta\cdot x}\theta_j$ for $j=1,2$.
\end{remark}

\subsection{Born approximation for fixed-angle data.}
The definition of the Born approximation in this case is more complicated than in the previous one. The reason is that in the backscattering case, since $\zeta=-\theta$, we can chose the frequencies of the incident waves in an appropriate way so that we can combine \eqref{amplitudppR}--\eqref{amplitudssR} to recover $\widehat{\mathrm{Q}}(-2 \omega \theta)e_i$ (see \eqref{gorroQ}, \eqref{b_p}, \eqref{b_s}, Definition \ref{BornB} and Remark \ref{datosB}). In fact, we have recovered $\widehat{\mathrm{Q}}(\xi)e_i$ for all $\xi\in\mathbb{R}^2$, since $\xi$ can be written uniquely as $\xi=-2\omega\theta$ with $\omega>0$ and $\theta\in\mathbb{S}^1$. However, in the fixed angle case, where $\theta$ is fixed, this is not possible because, there does not exist $\zeta\in\mathbb{S}^1$ such that 
$\xi=\omega_1(\zeta-\theta)=\omega_2(K\zeta-\theta)$ with $\omega_1,\omega_2>0$, except for the trivial case $K=1.$
But we do have that there exist two different $\zeta_1$ and $\zeta_2$ in $\mathbb{S}^1$ such that $\xi=\omega_1(\zeta_1-\theta)=\omega_2(K\zeta_2-\theta)$,
at least for the case $K>1$ whenever $\xi\cdot\theta<0$.
This will allow us to recover $\widehat{\mathrm{Q}}(\xi)e_i$ but in a different and more complicated way. 

In this case, given a fixed angle $\theta\in \mathbb{S}^1$, we look for a matrix depending on $\theta$, called the Born approximation for fixed-angle data, written in terms of the scattering amplitudes. As in the previous subsection, to do that it is enough to find expressions for the vectors $\mathrm{Q}(x)\theta$ and $\mathrm{Q}(x)\theta^\perp$.
Following \cite[Section 4]{bfprv1},
we will try to write the vector $\mathrm{Q}(x)\theta$ in terms of the scattering amplitudes, which are written in terms of the Fourier transform of $\mathrm{Q}$, and therefore, we start writing
\begin{equation*}
\label{Qtheta}
\mathrm{Q}(x)\theta
=
\int_{\mathbb{R}^2}e^{ix\cdot\xi}\widehat{\mathrm{Q}}(\xi)\theta\,d\xi.
\end{equation*}
Now, 
if we consider as incident p-wave
$\mathbf{u}_{i}^{p}(x)=e^{ik_{p}\theta\cdot x}\theta$
and as incident s-wave 
$\mathbf{u}_{i}^s(x)=e^{ik_{s}\theta^\perp\cdot x}\theta$, 
from 
\eqref{amplitudppR}--\eqref{amplitudssR} with $\varphi=\theta^\perp$ and $\phi=\theta$ we get
\begin{align}
\label{amplitudppA}
\Pi_{\zeta}\widehat{\mathrm{Q}}\left(\omega(\zeta-\theta)\right)\theta
=
&
(2\mu+\lambda)\mathbf{v}_{p,\infty}^{p}(\zeta;\sqrt{2\mu+\lambda}\omega,\theta)
\\
&
-
\Pi_{\zeta}(\mathrm{Q}\mathbf{v}^{p}(\cdot;\sqrt{2\mu+\lambda}\omega,\theta))\ \widehat{}\ \left(\omega\,\zeta\right),
\nonumber
\end{align}
\begin{align}
\label{amplitudpsA}
(\mathrm{I}-\Pi_{\zeta})\widehat{\mathrm{Q}}\left(\omega(K\zeta-\theta)\right)\theta
=
&
\mu\mathbf{v}_{s,\infty}^{p}(\zeta;\sqrt{2\mu+\lambda}\omega,\theta)
\\
&
-
(\mathrm{I}-\Pi_{\zeta})(\mathrm{Q}\mathbf{v}^{p}(\cdot;\sqrt{2\mu+\lambda}\omega,\theta))\ \widehat{}\ \left(K\omega\zeta\right),
\nonumber
\end{align}
\begin{align}
\label{amplitudspA}
\Pi_{\zeta}\widehat{\mathrm{Q}}\left(\omega(K^{-1}\zeta-\theta^\perp)\right)\theta
=
&
(2\mu+\lambda)\mathbf{v}_{p,\infty}^{s}(\zeta;\sqrt{\mu}\omega,\theta^\perp,\theta)
\\
&
-
\Pi_{\zeta}(\mathrm{Q}\mathbf{v}^{s}(\cdot;\sqrt{\mu}\omega,\theta^\perp,\theta))\ \widehat{}\ \left(K^{-1}\omega\,\zeta\right),
\nonumber
\end{align}
\begin{align}
\label{amplitudssA}
(\mathrm{I}-\Pi_{\zeta})\widehat{\mathrm{Q}}\left(\omega(\zeta-\theta^\perp)\right)\theta
=
&
\mu\mathbf{v}_{s,\infty}^{s}(\zeta;\sqrt{\mu}\omega,\theta^\perp,\theta)
\\
&
-
(\mathrm{I}-\Pi_{\zeta})(\mathrm{Q}\mathbf{v}^{s}(\cdot;\sqrt{\mu}\omega,\theta^\perp,\theta))\ \widehat{}\ \left(\omega\zeta\right).
\nonumber
\end{align}
Notice that for $K>1,$ every $\xi$ in the half plane 
$H_{\theta}=\{\xi\in\mathbb{R}^{2}:\xi\cdot\theta<0\}$
can be uniquely represented as 
\begin{equation}
\xi=\omega_{1}(\zeta_{1}-\theta)=\omega_{2}(K\zeta_{2}-\theta), 
\label{representaxi}%
\end{equation}
with $\omega_1,\omega_2>0$ and $\zeta_1,\zeta_2\in \mathbb{S}^1$ given by
\begin{equation}
\label{w1z1}
\omega_{1}
=
-\frac{|\xi|^2}{2\theta\cdot\xi}=\omega_{1}^\theta\,,
\quad
\zeta_{1}
=
-\frac{(2\xi\cdot\theta)}{|\xi|^2}\,\xi+\theta=\zeta_{1}^\theta\,
\end{equation}
\begin{equation}
\label{w2}
\omega_{2}
=
\frac{|\xi|^2}{-\xi\cdot\theta+\sqrt{(\xi\cdot\theta)^{2}+|\xi|^2(K^{2}-1)}}
=\omega_{2}^{\theta,K},
\end{equation}
and
\begin{equation}
\label{z2}
\zeta_{2}
=
\frac{-\xi\cdot\theta+\sqrt{\left(\xi\cdot\theta\right)^{2}+|\xi|^2(K^{2}-1)}}{|\xi|^2}\,
\frac{\xi}{K}+\frac{\theta}{K}
=\zeta_{2}^{\theta,K}.
\end{equation}
The sets $\displaystyle \left\{\xi\in\mathbb{R}^2\,/\,\xi=\omega_1(\zeta_1-\theta)\right\}$ and
$\displaystyle \left\{\xi\in\mathbb{R}^2\,/\,\xi=\omega_2(K\zeta_2-\theta)\right\}$ are known as Ewald and $K$-Ewald spheres (see Figure \ref{EsferasEwald} below). 
As can be seen in the first and second picture of Figure \ref{EsferasEwald}, 
any $\xi$ in the half plane $H_\theta$ belongs to a single Ewald sphere and to a single $K$-Ewald sphere whenever $K>1$.
But as the third picture in Figure \ref{EsferasEwald} shows,
this is not longer true in case $K<1$, and therefore, the definition of the Born approximation for fixed-angle data will be different depending on whether $K\ge 1$ or not. So we distinguish two cases: $K\ge 1$ and $K\le 1$.

\begin{figure}
$$
\begin{array}{c}
\begin{picture}(10,150)
\scalebox{.8}{
	\rput(0.5,-0.7){\normalsize{\textcolor[rgb]{0.68,0.13,0.04}{$\theta$}}}
	\rput(-2.5,5.5){\normalsize{$\theta,\zeta_1\in \mathbb{S}^1$}}
	\rput(-2.5,5){\normalsize{$\omega_1>0$}}
	\rput(-3,0.5){\normalsize{$\displaystyle\alpha= \frac{\pi}{4}$}}
	\rput(3,5.5){Ewald spheres}
	\rput(3,5){\normalsize{$\displaystyle \left\{\xi\in\mathbb{R}^2\,/\,\xi=\omega_1(\zeta_1-\theta)\right\}$}}
	\psaxes[linewidth=.7pt,labels=none,ticks=none]{->}(0,0)(-4.5,-1.5)(4.5,6)
	\psline[linecolor=rojo,linewidth=2pt]{->}(0,0)(0,-1)
	\psline[linecolor=gray,linewidth=1pt](0,0)(4,4)
	\psline[linecolor=gray,linewidth=1pt](0,0)(-4,4)
	\psarc[linecolor=gray](0,0){0.8}{90}{135}
	\psarc[linecolor=gray](0,0){0.87}{90}{135}
	\rput(-0.2,0.5){\normalsize{\textcolor[rgb]{0.63,0.63,0.63}{$\alpha$}}}
	\psline[linecolor=gray,linestyle=dashed,linewidth=1pt,dash=3pt 6pt](0,1)(1,1)
	\psline[linecolor=gray,linestyle=dashed,linewidth=1pt,dash=3pt 6pt](0,1.5)(1.5,1.5)
	\psline[linecolor=gray,linestyle=dashed,linewidth=1pt,dash=3pt 6pt](0,2)(2,2)
	\psline[linecolor=gray,linestyle=dashed,linewidth=1pt,dash=3pt 6pt](0,2.5)(2.5,2.5)
	\pscircle[linecolor=blue](0,1){1}
	\pscircle[linecolor=blue,fillstyle=solid,fillcolor=blue](0,1){0.07}
	\pscircle[linecolor=blue!75](0,1.5){1.5}
	\pscircle[linecolor=blue!75,fillstyle=solid,fillcolor=blue!75](0,1.5){0.07}
	\pscircle[linecolor=blue!50](0,2){2}
	\pscircle[linecolor=blue!50,fillstyle=solid,fillcolor=blue!50](0,2){0.07}
	\pscircle[linecolor=blue!25](0,2.5){2.5}
	\pscircle[linecolor=blue!25,fillstyle=solid,fillcolor=blue!25](0,2.5){0.07}
}
\label{pp}
\end{picture}
\\
\begin{picture}(10,155)
\scalebox{.8}{
	\rput(0.5,-0.7){\normalsize{\textcolor[rgb]{0.68,0.13,0.04}{$\theta$}}}
	\rput(-2.5,4.75){\normalsize{$\theta,\zeta_2\in \mathbb{S}^1$}}
	\rput(-2.5,4.25){\normalsize{$\omega_2>0$}}
	\rput(-4,0.5){\normalsize{$\displaystyle\alpha= \arctan K\in\left(\frac{\pi}{4},\frac{\pi}{2}\right)$}}
	\rput(3,4.75){$K-$Ewald spheres for $K>1$}
	\rput(3,4.25){\normalsize{$\displaystyle \left\{\xi\in\mathbb{R}^2\,/\,\xi=\omega_2(K\zeta_2-\theta)\right\}$}}
	\psaxes[linewidth=.7pt,labels=none,ticks=none]{->}(0,0)(-4.5,-2)(4.5,5)
	\psline[linecolor=rojo,linewidth=2pt]{->}(0,0)(0,-1)
	\psline[linecolor=gray,linewidth=1pt](0,0)(4,2)
	\psline[linecolor=gray,linewidth=1pt](0,0)(-4,2)
	\psarc[linecolor=gray](0,0){0.8}{90}{150}
	\psarc[linecolor=gray](0,0){0.87}{90}{150}
	\rput(-0.3,0.5){\normalsize{\textcolor[rgb]{0.63,0.63,0.63}{$\alpha$}}}
	\psline[linecolor=gray,linestyle=dashed,linewidth=1pt,dash=3pt 6pt](0,0.5)(1,0.5)
	\psline[linecolor=gray,linestyle=dashed,linewidth=1pt,dash=3pt 6pt](0,0.75)(1.5,0.75)
	\psline[linecolor=gray,linestyle=dashed,linewidth=1pt,dash=3pt 6pt](0,1)(2,1)
	\psline[linecolor=gray,linestyle=dashed,linewidth=1pt,dash=3pt 6pt](0,1.25)(2.5,1.25)
	\pscircle[linecolor=blue](0,0.5){1}
	\pscircle[linecolor=blue,fillstyle=solid,fillcolor=blue](0,0.5){0.07}
	\pscircle[linecolor=blue!75](0,0.75){1.5}
	\pscircle[linecolor=blue!75,fillstyle=solid,fillcolor=blue!75](0,0.75){0.07}
	\pscircle[linecolor=blue!50](0,1){2}
	\pscircle[linecolor=blue!50,fillstyle=solid,fillcolor=blue!50](0,1){0.07}
	\pscircle[linecolor=blue!25](0,1.25){2.5}
	\pscircle[linecolor=blue!25,fillstyle=solid,fillcolor=blue!25](0,1.25){0.07}
}
\end{picture}
\\
\begin{picture}(10,210)
\scalebox{.8}{
	\rput(0.5,-0.7){\normalsize{\textcolor[rgb]{0.68,0.13,0.04}{$\theta$}}}
	\rput(-3.5,2.5){\normalsize{$\theta,\zeta_2\in \mathbb{S}^1$}}
	\rput(-3.5,2){\normalsize{$\omega_2>0$}}
	\rput(-2.5,0.5){\normalsize{$\displaystyle\alpha= \arctan K^{-1}\in\left(0,\frac{\pi}{4}\right)$}}
	\rput(4,2.5){$K-$Ewald spheres for $K<1$}
	\rput(4,2){\normalsize{$\displaystyle \left\{\xi\in\mathbb{R}^2\,/\,\xi=\omega_2(K\zeta_2-\theta)\right\}$}}
	\psaxes[linewidth=.7pt,labels=none,ticks=none]{->}(0,0)(-4.5,-1.5)(4.5,7)
	\psline[linecolor=rojo,linewidth=2pt]{->}(0,0)(0,-1)
	\psline[linecolor=gray,linewidth=1pt](0,0)(3,6)
	\psline[linecolor=gray,linewidth=1pt](0,0)(-3,6)
	\psarc[linecolor=gray](0,0){0.9}{90}{120}
	\psarc[linecolor=gray](0,0){0.97}{90}{120}
	\rput(-0.15,0.7){\normalsize{\textcolor[rgb]{0.63,0.63,0.63}{$\alpha$}}}
	\psline[linecolor=gray,linestyle=dashed,linewidth=1pt,dash=3pt 6pt](0,1)(0.5,1)
	\psline[linecolor=gray,linestyle=dashed,linewidth=1pt,dash=3pt 6pt](0,2)(1,2)
	\psline[linecolor=gray,linestyle=dashed,linewidth=1pt,dash=3pt 6pt](0,3)(1.5,3)
	\psline[linecolor=gray,linestyle=dashed,linewidth=1pt,dash=3pt 6pt](0,4)(2,4)
	\pscircle[linecolor=blue](0,1){0.5}
	\pscircle[linecolor=blue,fillstyle=solid,fillcolor=blue](0,1){0.07}
	\pscircle[linecolor=blue!75](0,2){1}
	\pscircle[linecolor=blue!75,fillstyle=solid,fillcolor=blue!75](0,2){0.07}
	\pscircle[linecolor=blue!50](0,3){1.5}
	\pscircle[linecolor=blue!50,fillstyle=solid,fillcolor=blue!50](0,3){0.07}
	\pscircle[linecolor=blue!25](0,4){2}
	\pscircle[linecolor=blue!25,fillstyle=solid,fillcolor=blue!25](0,4){0.07}
}
\end{picture}
\end{array}
$$
\vspace*{0.5cm}
	\caption{
	$K-$Ewald spheres for different values of $K=\sqrt{2\mu+\lambda}/\sqrt{\mu}$.
}
	\label{EsferasEwald}
\end{figure}
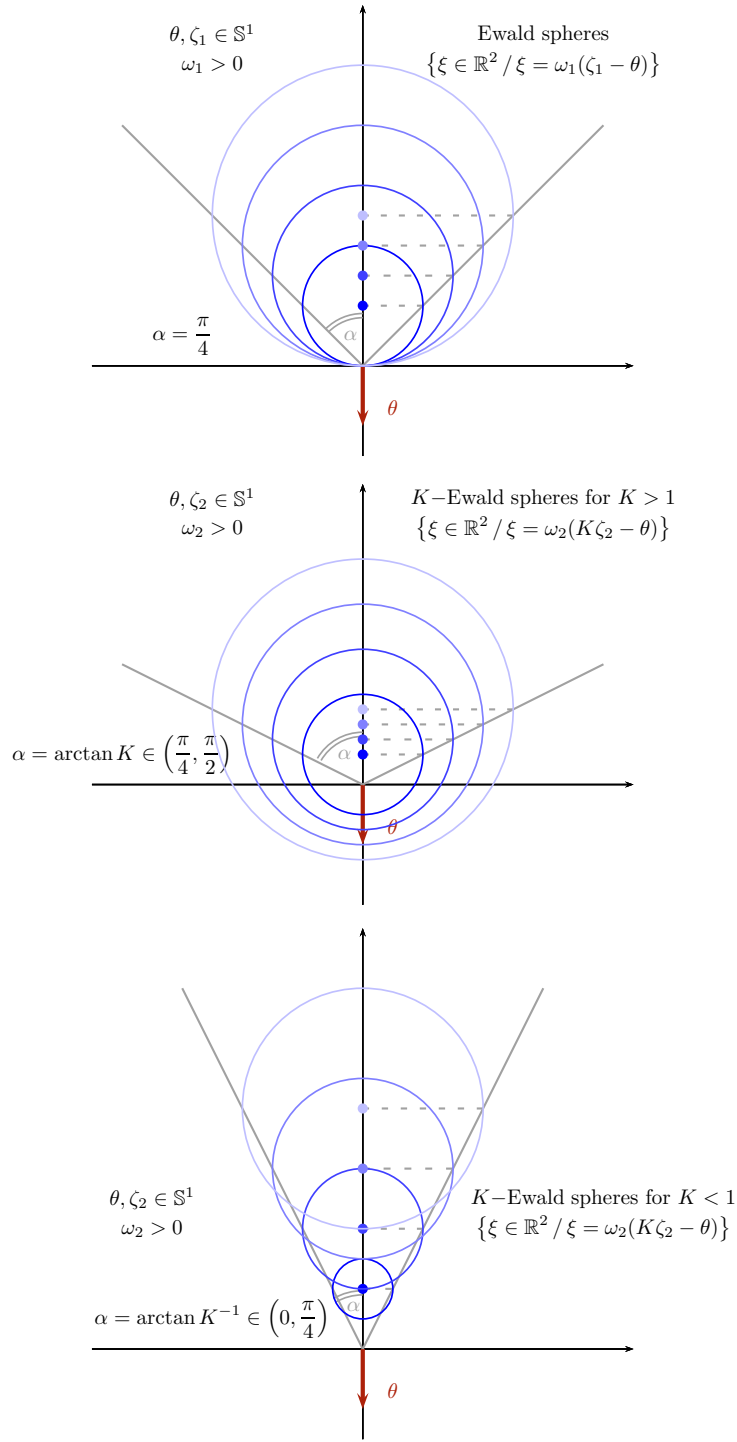

\underline{Case $K>1$}.
As we said before, for $K> 1$ and $\xi\in H_\theta$, identities \eqref{representaxi}--\eqref{z2} hold. Moreover, we have that $\zeta_{1}\cdot\zeta_{2}\geq 1/K\neq 0$, and thus any vector $\mathbf{v}$ can be written in terms of the orthogonal projections $\Pi_{\zeta_1}$ and $\Pi_{\zeta_2}$ as follows
$$
\mathbf{v}
=
(\mathrm{I}-\Pi_{\zeta_{2}})\mathbf{v}
+
\left[
\left(\Pi_{\zeta_{1}}\mathbf{v}-(\mathrm{I}-\Pi_{\zeta_{2}})\mathbf{v}\right)
\cdot\frac{\zeta_{1}}{\zeta_{1}\cdot\zeta_{2}}
\right]
\zeta_{2}.
$$
Applying this expression to the vector $\mathbf{v}=\widehat{\mathrm{Q}}(\xi)\theta$ with $\xi\in H_\theta$, and using \eqref{amplitudppA} with $\omega=\omega_1$ and $\zeta=\zeta_1$, and \eqref{amplitudpsA} with $\omega=\omega_2$ and $\zeta=\zeta_1$, we have that
$$
\widehat{\mathrm{Q}}(\xi)\,\theta=\mathbf{v}_{\infty}^{\theta,p}(\xi)-
\mathbf{e}^{\theta,p}(\xi),
$$
where $\mathbf{v}_{\infty}^{\theta,p}$ is written in terms of the $p\rightarrow p$ and the $p\rightarrow s$ scattering amplitudes as follows,
\begin{multline}
\label{amplitudp}
\mathbf{v}_{\infty}^{\theta,p}(\xi)
=
\mu\mathbf{v}_{s,\infty}^{p}(\zeta_2;\sqrt{2\mu+\lambda}\omega_2,\theta)
\\
+
\left(
(2\mu+\lambda)\mathbf{v}_{p,\infty}^{p}(\zeta_1;\sqrt{2\mu+\lambda}\omega_1,\theta)
-
\mu\mathbf{v}_{s,\infty}^{p}(\zeta_2;\sqrt{2\mu+\lambda}\omega_2,\theta)
\right)
\cdot\frac{\zeta_{1}}{\zeta_{1}\cdot\zeta_{2}}
\,\zeta_{2},
\end{multline}
with $\omega_1=\omega_1^{\theta}(\xi)$ and $\zeta_1=\zeta_1^{\theta}(\xi)$ given in \eqref{w1z1} and,
$\omega_2=\omega_2^{\theta,K}(\xi)$ and $\zeta_2=\zeta_2^{\theta,K}(\xi)$ given in \eqref{w2} and \eqref{z2} respectively,
and $\mathbf{e}^{\theta,p}$ is an error term given by
\begin{multline}
\label{errorp}
\mathbf{e}^{\theta,p}(\xi)
=
(\mathrm{I}-\Pi_{\zeta_{2}})(\mathrm{Q}\mathbf{v}^{p}(\cdot;\sqrt{2\mu+\lambda}\omega_2,\theta))\ \widehat{}\ (K\omega_2\,\zeta_2)
\\
+
\left[
\left(
\Pi_{\zeta_{1}}(\mathrm{Q}\mathbf{v}^{p}(\cdot;\sqrt{2\mu+\lambda}\omega_1,\theta))\ \widehat{}\ (\omega_1\,\zeta_1)
\right.\right.
\\
\left.\left.
-
(\mathrm{I}-\Pi_{\zeta_{2}})(\mathrm{Q}\mathbf{v}^{p}(\cdot;\sqrt{2\mu+\lambda}\omega_2,\theta))\ \widehat{}\ (K\omega_2\,\zeta_2)
\right)
\cdot\frac{\zeta_{1}}{\zeta_{1}\cdot\zeta_{2}}
\right]
\zeta_{2}.
\end{multline}

In a similar way, we can write the vector $\mathbf{v}=\widehat{\mathrm{Q}}(\xi)\theta$ with $\xi\in H_{-\theta}$ in the form
$$
\widehat{\mathrm{Q}}(\xi)\,\theta=-\mathbf{v}_{\infty}^{-\theta,p}(\xi)+\mathbf{e}^{-\theta,p}(\xi),
$$
and therefore, since $\mathbb{R}^2\setminus\{\xi\in\mathbb{R}^2/\xi\cdot\theta=0\}=H_\theta\dot\cup H_{-\theta}$, whenever $\xi\cdot\theta\neq 0$ we have that
\begin{equation}
\label{Qgorro_theta}
\widehat{\mathrm{Q}}(\xi)\,\theta
=
\left(\mathbf{v}_{\infty}^{\theta,p}(\xi)-\mathbf{e}^{\theta,p}(\xi)\right)\chi_{H_\theta}(\xi)
-
\left(\mathbf{v}_{\infty}^{-\theta,p}(\xi)-\mathbf{e}^{-\theta,p}(\xi)\right)\chi_{H_{-\theta}}(\xi).
\end{equation}
Arguing as before but replacing $\theta$ by $\theta^\perp$, whenever $\xi\cdot\theta^\perp\neq 0$ we get
\begin{equation}
\label{Qgorro_theta_orto}
\widehat{\mathrm{Q}}(\xi)\,\theta^\perp
=
\left(\mathbf{v}_{\infty}^{\theta^\perp,p}(\xi)-\mathbf{e}^{\theta^\perp,p}(\xi)\right)
\chi_{H_{\theta^\perp}}(\xi)
-
\left(\mathbf{v}_{\infty}^{-\theta^\perp,p}(\xi)-\mathbf{e}^{-\theta^\perp,p}(\xi)\right)
\chi_{H_{-\theta^\perp}}(\xi).
\end{equation}
Thus, in the case $K\ge 1$, it is natural to define the Born approximation for fixed-angle data as follows:
\begin{definition}
	\label{BornAK>1}
	Let $K>1$.
	Given $\mathrm{Q}$ a matrix of order $2$ and an angle $\theta\in\mathbb{S}^1$,
	we define the \emph{Born approximation for fixed-angle data} of $\mathrm{Q}$ as a matrix $\mathrm{Q}_B^\theta$ such that for any $\xi\in\mathbb{R}^2\setminus\{\xi\in\mathbb{R}^2/\xi\cdot\theta=0 \textrm{ or } \xi\cdot\theta^\perp=0\}$,
	\begin{multline}
	\widehat{\mathrm{\mathrm{Q}}_{B}^\theta}(\xi)e_{i}
	=
	( e_{i}\cdot \theta)
	\left(
	\mathbf{v}_{\infty}^{\theta,p}(\xi)\chi_{H_\theta}(\xi)
	-
	\mathbf{v}_{\infty}^{-\theta,p}(\xi)\chi_{H_{-\theta}}(\xi)
	\right)
	\\
	+
	( e_{i}\cdot \theta^\perp)
	\left(
	\mathbf{v}_{\infty}^{\theta^\perp,p}(\xi)\chi_{H_{\theta^\perp}}(\xi)
	-
	\mathbf{v}_{\infty}^{-\theta^\perp,p}(\xi)\chi_{H_{-\theta^\perp}}(\xi)
	\right),
	\
	i=1,2,
	\label{Def_BornAK>1}
	\end{multline}
	 where $\{e_1,  e_2\}$ is the canonical base of $\mathbb{R}^{2},$ and 
	$\mathbf{v}_{\infty}^{\theta,p}$ is given in \eqref{amplitudp}.
\end{definition}
\begin{remark}
	\label{rectasortogonales}
	Note that although in the previous definition $\xi\cdot\theta=0$ and $\xi\cdot\theta^\perp=0$ are not allowed,
	$\mathrm{Q}_B^\theta$ is defined for all $x\in\mathbb{R}^2$, because
	if we write $X_\theta=\mathbb{R}^2\setminus\{\xi\in\mathbb{R}^2/\xi\cdot\theta=0 \textrm{ or } \xi\cdot\theta^\perp=0\}$ we have that
	\begin{equation*}
	\mathrm{Q}_B^\theta(x)e_i 
	=\int_{\mathbb{R}^2}
	e^{i \xi \cdot x}\widehat{\mathrm{Q}_B^\theta}(\xi)e_i\,d\xi
	=\int_{X_\theta} 
	e^{i \xi \cdot x}\widehat{\mathrm{Q}_B^\theta}(\xi)e_i\,d\xi.
	\end{equation*} 
\end{remark}
From the Definition \ref{BornAK>1} and the identities \eqref{Qgorro_theta} and \eqref{Qgorro_theta_orto} we have that 
\begin{equation}
\label{Q_QB_E}
\mathrm{Q}(x)=\mathrm{Q}_B^\theta(x)-\mathrm{E}^\theta(x),\qquad x\in\mathbb{R}^2,
\end{equation}
where $\mathrm{E}^\theta=\mathrm{E}^\theta({\mathrm{Q}})$ is the error term defined by  	
\begin{multline}
\widehat{\mathrm{\mathrm{E}}^\theta}(\xi)e_{i}
=
( e_{i}\cdot \theta)
\left(
\mathbf{e}^{\theta,p}(\xi)\chi_{H_\theta}(\xi)-\mathbf{e}^{-\theta,p}(\xi)\chi_{H_{-\theta}}(\xi)
\right)
\\
+
( e_{i}\cdot \theta^\perp)
\left(
\mathbf{e}^{\theta^\perp,p}(\xi)\chi_{H_{\theta^\perp}}(\xi)-\mathbf{e}^{-\theta^\perp,p}(\xi)\chi_{H_{-\theta^\perp}}(\xi)
\right),
\
i=1,2,
\label{errorAK>1}
\end{multline}
with $\xi\in\mathbb{R}^2\setminus\{\xi\in\mathbb{R}^2/\xi\cdot\theta=0 \textrm{ or } \xi\cdot\theta^\perp=0\}$ and
$\mathbf{e}^{\theta,p}$ given in \eqref{errorp}.

\begin{remark}
	\label{dimensionAF}
This definition can be extended to dimension $n>2$  (see \cite[Definition 4.1]{bfprv1}), and notice that it depends on the choice of $\theta^\perp$.
\end{remark}

\underline{Case $K<1$}.
If $K<1$ we cannot use any more longitudinal incident waves, because in this case
we cannot find $\omega_2>0$ and $\zeta_2\in \mathbb{S}^1$ such that any $\xi\in H_{\theta}$
can be represented as
$\xi=\omega_{2}(K\zeta_{2}-\theta)$.
This is equivalent to saying that we cannot cover the half-plane $H_{\theta}$ with $K$-Ewald spheres when $K<1$, as shown in the third picture in Figure \ref{EsferasEwald}.
But we can use transverse incident waves to define the Born approximation. 

More precisely, we will use \eqref{amplitudspA} and \eqref{amplitudssA}, since in this case, every $\xi\in H_{\theta^\perp}$ can be uniquely represented as 
\begin{equation}
\xi=\tilde{\omega}_{1}(\tilde{\zeta}_{1}-\theta^\perp)=
\tilde{\omega}_{2}(K^{-1}\tilde{\zeta}_{2}-\theta^\perp), 
\label{representaxiK<1}%
\end{equation}
with $\tilde{\omega}_1=\omega_1^{\theta^\perp}(\xi)$ and 
$\tilde{\zeta}_1=\zeta_1^{\theta^\perp}(\xi)$ given by \eqref{w1z1} and,
$\tilde{\omega}_2=\omega_2^{\theta^\perp,K^{-1}}(\xi)$ and $\tilde{\zeta}_2=\zeta_2^{\theta^\perp,K^{-1}}(\xi)$ given by \eqref{w2} and \eqref{z2} respectively.
And therefore, arguing as in the other case,
for any $\xi\in\mathbb{R}^2\setminus\{\xi\in\mathbb{R}^2/\xi\cdot\theta=0 \textrm{ or } \xi\cdot\theta^\perp=0\}$ we have that
\begin{equation}
\label{Qgorro_thetaK<1}
\widehat{\mathrm{Q}}(\xi)\,\theta
=
\left(\mathbf{v}_{\infty}^{\theta^\perp,s}(\xi)-\mathbf{e}^{\theta^\perp,s}(\xi)\right)\chi_{H_{\theta^\perp}}(\xi)
-
\left(\mathbf{v}_{\infty}^{-\theta^\perp,s}(\xi)-\mathbf{e}^{-\theta^\perp,s}(\xi)\right)\chi_{H_{-\theta^\perp}}(\xi),
\end{equation}
and
\begin{equation}
\label{Qgorro_theta_ortoK<1}
\widehat{\mathrm{Q}}(\xi)\,\theta^\perp
=
\left(\mathbf{v}_{\infty}^{-\theta,s}(\xi)-\mathbf{e}^{-\theta,s}(\xi)\right)
\chi_{H_{-\theta}}(\xi)
-
\left(\mathbf{v}_{\infty}^{\theta,s}(\xi)-\mathbf{e}^{\theta,s}(\xi)\right)
\chi_{H_{\theta}}(\xi),
\end{equation}
where
\begin{multline}
\label{amplituds}
\mathbf{v}_{\infty}^{\theta^\perp, s}(\xi)
=
\mu\mathbf{v}_{s,\infty}^{s}(\tilde{\zeta}_1;\sqrt{\mu}\tilde{\omega}_1,\theta^\perp,\theta)
\\
+
\left(
(2\mu+\lambda)\mathbf{v}_{p,\infty}^{s}(\tilde{\zeta}_2;\sqrt{\mu}\tilde{\omega}_2,\theta^\perp,\theta)
-
\mu\mathbf{v}_{s,\infty}^{s}(\tilde{\zeta}_1;\sqrt{\mu}\tilde{\omega}_1,\theta^\perp,\theta)
\right)
\cdot\frac{\tilde{\zeta}_{2}}{\tilde{\zeta}_{1}\cdot\tilde{\zeta}_{2}}
\,\tilde{\zeta}_{1},
\end{multline}
and
\begin{multline}
\label{errors}
\mathbf{e}^{\theta^\perp,s}(\xi)
=
(\mathrm{I}-\Pi_{\tilde{\zeta}_{1}})
(\mathrm{Q}\mathbf{v}^{s}(\cdot;\sqrt{\mu}\tilde{\omega}_1,\theta^\perp,\theta))\ \widehat{}\ (\tilde{\omega}_1\,\tilde{\zeta}_1)
\\
+
\left[
\left(
\Pi_{\tilde{\zeta}_{2}}(\mathrm{Q}\mathbf{v}^{s}(\cdot;\sqrt{\mu}\tilde{\omega}_2,\theta^\perp,\theta))\ \widehat{}\ (K^{-1}\tilde{\omega}_2\,\tilde{\zeta}_2)
\right.\right.
\\
\left.\left.
-
(\mathrm{I}-\Pi_{\tilde{\zeta}_{1}})
(\mathrm{Q}\mathbf{v}^{s}(\cdot;\sqrt{\mu}\tilde{\omega}_1,\theta^\perp,\theta))\ \widehat{}\ (\tilde{\omega}_1\,\tilde{\zeta}_1)
\right)
\cdot\frac{\tilde{\zeta}_{2}}{\tilde{\zeta}_{1}\cdot\tilde{\zeta}_{2}}
\right]
\tilde{\zeta}_{1}.
\end{multline}
Thus the definition in this case is as follows.
\begin{definition}
	\label{BornAK<1}
	Let $K<1$.
	Given $\mathrm{Q}$ a matrix of order $2$ and an angle $\theta\in\mathbb{S}^1$,
	we define the \emph{Born approximation for fixed-angle data} of $\mathrm{Q}$ as a matrix $\mathrm{Q}_B^\theta$ such that for any $\xi\in\mathbb{R}^2\setminus\{\xi\in\mathbb{R}^2/\xi\cdot\theta=0 \textrm{ or } \xi\cdot\theta^\perp=0\}\footnote{See Remark \ref{rectasortogonales}.}$,
	\begin{multline}
	\widehat{\mathrm{\mathrm{Q}}_{B}^\theta}(\xi)e_{i}
	=
	( e_{i}\cdot \theta)
	\left(
	\mathbf{v}_{\infty}^{\theta^\perp,s}(\xi)\chi_{H_{\theta^\perp}}(\xi)
	-
	\mathbf{v}_{\infty}^{-\theta^\perp,s}(\xi)\chi_{H_{-\theta^\perp}}(\xi)
	\right)
	\\
	+
	( e_{i}\cdot \theta^\perp)
	\left(
	\mathbf{v}_{\infty}^{-\theta,s}(\xi)\chi_{H_{-\theta}}(\xi)
	-
	\mathbf{v}_{\infty}^{\theta,s}(\xi)\chi_{H_{\theta}}(\xi)
	\right),
	\
	i=1,2,
	\label{Def_BornAK<1}
	\end{multline}
	where $\{e_1,  e_2\}$ is the canonical base of $\mathbb{R}^{2},$ and 
	$\mathbf{v}_{\infty}^{\theta^\perp,s}$ is given in \eqref{amplituds}.
\end{definition}
From this definition, \eqref{Qgorro_thetaK<1} and \eqref{Qgorro_theta_ortoK<1} we have that 
\begin{equation}
\label{Q_QB_EK<1}
\mathrm{Q}(x)=\mathrm{Q}_B^\theta(x)-\mathrm{E}^\theta(x),\qquad x\in\mathbb{R}^2,
\end{equation}
where $\mathrm{E}^\theta=\mathrm{E}^\theta({\mathrm{Q}})$ is the error term defined by 	
\begin{multline}
\widehat{\mathrm{\mathrm{E}}^\theta}(\xi)e_{i}
=
( e_{i}\cdot \theta)
\left(
\mathbf{e}^{\theta^\perp,s}(\xi)\chi_{H_{\theta^\perp}}(\xi)
-
\mathbf{e}^{-\theta^\perp,s}(\xi)\chi_{H_{-\theta^\perp}}(\xi)
\right)
\\
+
( e_{i}\cdot \theta^\perp)
\left(
\mathbf{e}^{-\theta,s}(\xi)\chi_{H_{-\theta}}(\xi)
-
\mathbf{e}^{\theta,s}(\xi)\chi_{H_{\theta}}(\xi)
\right),
\
i=1,2,
\label{errorAK<1}
\end{multline}
with $\xi\in\mathbb{R}^2\setminus\{\xi\in\mathbb{R}^2/\xi\cdot\theta=0 \textrm{ or } \xi\cdot\theta^\perp=0\}$ 
and
$\mathbf{e}^{\theta^\perp,s}$ given in \eqref{errors}.
\begin{remark}
	\label{data}
	We note that if $k_s>k_p$ only $p\rightarrow p$ and $p\rightarrow s$ data are needed to define the Born approximation, while for $k_s<k_p$ only $s\rightarrow p$ and $s\rightarrow s$ data are needed. 
	
	Moreover, if $k_s=k_p$ we can choose and use only the $p\rightarrow p$ and $p\rightarrow s$ data to define the Born approximation, or use only the $s\rightarrow p$ and $s\rightarrow s$ data. Thus in this case, we have two different Born approximations. 
	This does not happen in the backscattering problem where all the scattering data are needed (see Definition \ref{BornB}).
	
	We would like to point out that both definitions need four scattering data, but in the fixed-angle data case all of them are longitudinal or transverse scattering data depending on $K$.
\end{remark}
\section{Iterative algorithms}
\label{SeccionAlgoritmos}
In this section we explain the algorithms that allow us to recover the unknown load $\mathrm{Q}$ from its scattering data. 
They are iterative algorithms that construct a sequence of approximations where the first term is given by the corresponding Born approximation.

We split this section into two subsections. In the first one the known scattering data are the backscattering data, and in the second one are the fixed-angle scattering data. 

\subsection{Iterative algorithm for backscattering data.}
In this subsections we assume that $\mathrm{Q}$ is an unknown matrix but its Born approximation for backscattering data $\mathrm{Q}_B$ is known.

In \eqref{definicion_juan} we have an identity that writes the unknown load $\mathrm{Q}$ as the difference of $\mathrm{Q}_B$ and an error term named $\mathrm{E}$.
But this error term depends on the unknown load itself. 
Actually $\mathrm{E}$ is defined by \eqref{errorb}, where it is written in terms of the errors $\mathbf{e}^p$ and $\mathbf{e}^s$ given in \eqref{hp} and \eqref{hs} respectively. And since these errors are written in terms of $\mathrm{Q}$, we have that
$\mathrm{E}=\mathrm{E}(\mathrm{Q})$. 

The idea of the iterative algorithm is to approximate this error term $\mathrm{E}(\mathrm{Q})$ by replacing the unknown matrix $\mathrm{Q}$ with a good approximation of it. As a first good approximation we have the Born approximation $\mathrm{Q}_B$, and in each iteration, we obtain a better approximation $\mathrm{Q}_n$ that we can use to approximate the error term of the next iteration. 

On the other hand, we have assumed that $\mathrm{Q}$ is compactly supported with support in $B(0,R)$ where $R>0$, but its Born approximation $\mathrm{Q}_B$, defined by \eqref{Def_BornB}, has not compact support. 
With this in mind, the iterative algorithm for the backscattering data is as follows:
\begin{equation}
\label{algoritmo_back}
\begin{aligned}
\mathrm{Q}_1&=\chi_{R} \,\mathrm{Q}_B
\\
\mathrm{Q}_{n+1}&=\chi_{R}\,\mathrm{Q}_B-\chi_R,\mathrm{E}_{n},
\qquad\qquad n=1,2,\ldots
\end{aligned}
\end{equation}
where $\mathrm{E}_{n}=\mathrm{E}(\mathrm{Q}_n)$ with $\mathrm{E}$ given by \eqref{errorb}.

\begin{figure}[t]
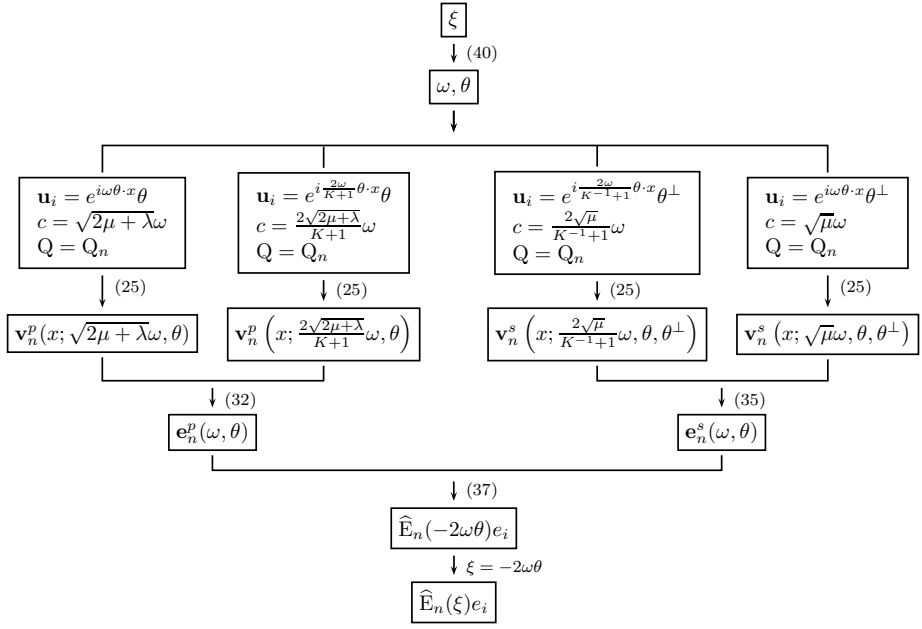

\centerline{
\scalebox{0.8}{
	\begin{minipage}{15cm}
$$
\begin{psmatrix}[colsep=-0.5cm,rowsep=0.5cm]
&&&\psframebox{\xi}&&&
\\
&&&\psframebox{\omega,\theta}&&&
\\[-1ex]
&&&{}&&&
\\
\psframebox{
\begin{array}{l}
\mathbf{u}_i=e^{i\omega\theta\cdot x}\theta\\
c=\sqrt{2\mu+\lambda}\omega\\
\mathrm{Q}=\mathrm{Q}_n
\end{array}}&&
\psframebox{
\begin{array}{l}
\mathbf{u}_i=e^{i\frac{2\omega}{K+1}\theta\cdot x}\theta\\
c=\frac{2\sqrt{2\mu+\lambda}}{K+1}\omega\\
\mathrm{Q}=\mathrm{Q}_n
\end{array}}&{}&
\psframebox{
\begin{array}{l}
\mathbf{u}_i=e^{i\frac{2\omega}{K^{-1}+1}\theta\cdot x}\theta^\perp\\
c=\frac{2\sqrt{\mu}}{K^{-1}+1}\omega\\
\mathrm{Q}=\mathrm{Q}_n
\end{array}}&&
\psframebox{
\begin{array}{l}
\mathbf{u}_i=e^{i\omega\theta\cdot x}\theta^\perp\\
c=\sqrt{\mu}\omega\\
\mathrm{Q}=\mathrm{Q}_n
\end{array}}
\\
\psframebox{\mathbf{v}^p_n(x;\sqrt{2\mu+\lambda}\omega,\theta)}
&&
\psframebox{\mathbf{v}^p_n\left(x;\frac{2\sqrt{2\mu+\lambda}}{K+1}\omega,\theta\right)}
&&
\psframebox{\mathbf{v}^s_n\left(x;\frac{2\sqrt{\mu}}{K^{-1}+1}\omega,\theta,\theta^\perp\right)}
&&
\psframebox{\mathbf{v}^s_n\left(x;\sqrt{\mu}\omega,\theta,\theta^\perp\right)}
\\[-2ex]
&{}&&{}&&{}&
\\
&\psframebox{\mathbf{e}^p_n(\omega,\theta)}&&{}&&
\psframebox{\mathbf{e}^s_n(\omega,\theta)}&
\\[-2ex]
&{}&&{}&&{}&
\\
&&&\psframebox{\widehat{\mathrm{E}}_n(-2\omega\theta)e_i}&&&
\\
&&&\psframebox{\widehat{\mathrm{E}}_n(\xi)e_i}&&&
\psset{nodesep=3pt}
\ncline{->}{1,4}{2,4}{\footnotesize\trput{\raisebox{0.0cm}{\eqref{omegatheta}}}}
\ncline{->}{2,4}{3,4}
\ncbar[angleA=90,nodesep=2pt,armA=0.45]{-}{4,1}{4,7}
\ncbar[angleA=90,nodesep=2pt,armA=0]{-}{4,3}{4,3}
\ncbar[angleA=90,nodesep=2pt,armA=0]{-}{4,5}{4,5}
\ncline{->}{2,4}{3,4}
\ncline{->}{4,1}{5,1}{\footnotesize\trput{\raisebox{-0.15cm}{\eqref{LSeq}}}}
\ncline{->}{4,3}{5,3}{\footnotesize\trput{\raisebox{-0.2cm}{\eqref{LSeq}}}}
\ncline{->}{4,5}{5,5}{\footnotesize\trput{\raisebox{-0.2cm}{\eqref{LSeq}}}}
\ncline{->}{4,7}{5,7}{\footnotesize\trput{\raisebox{-0.15cm}{\eqref{LSeq}}}}
\ncbar[angleA=-90,nodesep=2pt,armA=0.25,armB=0.25cm
]{-}{5,1}{5,3}
\ncbar[angleA=-90,nodesep=2pt,armA=0.25,armB=0.25cm
]{-}{5,5}{5,7}
\ncline{->}{6,2}{7,2}{\footnotesize\trput{\raisebox{0.05cm}{\eqref{hp}}}}
\ncline{->}{6,6}{7,6}{\footnotesize\trput{\raisebox{0.05cm}{\eqref{hs}}}}
\ncbar[angleA=-90,nodesep=2pt,armA=0.25,armB=0.25cm
]{-}{7,2}{7,6}
\ncline{->}{8,4}{9,4}{\footnotesize\trput{\raisebox{0.15cm}{\eqref{errorb}}}}
\ncline{->}{9,4}{10,4}{\footnotesize\trput{\xi=-2\omega\theta}}
\end{psmatrix}
$$ 
\end{minipage}
}
}
\caption{\label{diagramaEn}Construction of the Fourier transform of the error term at the n-th iteration for backscattering data.}
\end{figure}

The key point in each iteration is the construction of the error term. The diagram in Figure \ref{diagramaEn} shows the construction of the error at the $n$-th iteration.

In each iteration, 
given $\xi\in\mathbb{R}^2$,
we use \eqref{omegatheta} to get $\omega=\omega(\xi)$ and $\theta=\theta(\xi)$.
From here, we construct the four incident waves $\mathbf{u}_i$, two of them are plane p-waves and the other two are plane s-waves, and also the corresponding energies $c^2$. 
Then we introduce this information in the Lippmann-Schwinger equation given in \eqref{LSeq}, taking $\mathrm{Q}=\mathrm{Q}_n$. 
Solving the four equations we have obtained in this way we get the corresponding scattered solutions $\mathbf{v}_n^p$ or  $\mathbf{v}_n^s$. In order to solve these equations we use the algorithm described in \cite{BC}.
Once we have the scattered solutions we use the two that come from the plane p-waves in \eqref{hp} to construct the partial error term $\mathbf{e}_n^p$, and the two that come from the plane s-waves in \eqref{hs} to construct $\mathbf{e}_n^p$. In both identities, \eqref{hp} and \eqref{hs}, we take $\mathrm{Q}=\mathrm{Q}_n$.
Using these partial error terms in \eqref{errorb} we get the Fourier transform of the error at the n-th iteration  $\widehat{\mathrm{E}}_n(-2\omega\theta)e_i$.
This is just $\widehat{\mathrm{E}}_n(\xi)e_i$, since $\xi=-2\omega\theta$ (see Remark \ref{OmegaTheta}). 
Finally, from here we obtain $\mathrm{E}_n(x)e_i$ by taking the inverse Fourier transform.

\subsection{Iterative algorithm for fixed-angle data.}
In this subsection we again assume that
$\mathrm{Q}$ is unknown but now we know its Born approximation for fixed-angle data $\mathrm{Q}_B^\theta$.

We can argue as before because we have $\mathrm{Q}$ written as the difference between $\mathrm{Q}_B^\theta$ and an error term named $\mathrm{E}^\theta=\mathrm{E}^\theta(\mathrm{Q})$
in the identity \eqref{Q_QB_E} whenever $K>1$, and in the identity \eqref{Q_QB_EK<1} if $K<1$. And therefore, the iterative algorithm for the fixed-angle data is as follows:
\begin{equation}
\label{algoritmo_f_a}
\begin{aligned}
\mathrm{Q}_1^\theta&=\chi\,\mathrm{Q}_B^\theta
\\
\mathrm{Q}_{n+1}^\theta&=\chi\,\mathrm{Q}_B^\theta-\chi\,\mathrm{E}_{n}^\theta,
\qquad\qquad n=1,2,\ldots
\end{aligned}
\end{equation}
where $\chi$ is the cut-off function defined in the previous subsection and $\mathrm{E}_{n}^\theta=\mathrm{E}^\theta(\mathrm{Q}_n^\theta)$ with $\mathrm{E}^\theta$ given by \eqref{errorAK>1} if $K>1$ and by \eqref{errorAK<1} if $K<1$. 

As before, the key point in each iteration is the construction of the error term, but now this is more complicated.
The pictures in Figures \ref{cuadrantes} and \ref{diagrama_e_theta_p_n} show the construction of the Fourier transform of the error term at the $n$-th iteration for $K>1$. Similar diagrams can be constructed for the case $K<1$.

For $\theta$ fixed in $\mathbb{S}^1$, the lines $\xi\cdot\theta=0$ and $\xi\cdot\theta^\perp=0$ divide the space $\mathbb{R}^2$ into four quadrants: 
$H_\theta\cap H_{\theta^\perp}$, 
$H_{\theta}\cap H_{-\theta^\perp}$, 
$H_{-\theta}\cap H_{-\theta^\perp}$ and 
$H_{-\theta}\cap H_{\theta^\perp}$.
And given $\xi\in\mathbb{R}^2\setminus\{\xi\in\mathbb{R}^2/\xi\cdot\theta=0 \textrm{ or } \xi\cdot\theta^\perp=0\}$, there exists a unique quadrant where $\xi$ is.
As one can see in \eqref{errorAK>1} and \eqref{errorAK<1}, the construction of the error term is different depending on which quadrant $\xi$ is in.

More precisely, in the case $K>1$, if $\xi\in H_\theta\cap H_{\theta^\perp}$, from \eqref{errorAK>1} we have that just the partial error terms  
$\mathbf{e}_n^{\theta,p}(\xi)$ and $\mathbf{e}_n^{\theta^\perp,p}(\xi)$
are necessary to construct  
$\widehat{\mathrm{E}^\theta_n}(\xi)e_i$, and similarly if $\xi$ is in one of the other three quadrants (see Figure \ref{cuadrantes}).
Therefore, it is sufficient to see how to construct each of these partial error terms. 
Following the diagram in Figure \ref{diagrama_e_theta_p_n}, to get $\mathbf{e}_n^{\theta,p}(\xi)$, first of all we use 
\eqref{w1z1}, \eqref{w2} and \eqref{z2} to define 
$\omega_1^\theta$, $\zeta_1^\theta$, $\omega_2^{\theta,K}$ and $\zeta_2^{\theta,K}$. This can be done because 
$\xi\in H_\theta$. 
From here, we construct two incident waves $\mathbf{u}_i$ and the corresponding energies $c^2$.
Both waves are plane p-waves. 
Introducing this information in the Lippmann-Schwinger equation given in \eqref{LSeq} with $\mathrm{Q}=\mathrm{Q}_n^\theta$, and 
solving the two equations we obtain the two corresponding scattered solutions $\mathbf{v}_n^p$. Finally, using them in \eqref{errorp} we get $\mathbf{e}_n^{\theta,p}(\xi)$.
In a similar way we can obtain the other three partial errors, 
$\mathbf{e}_n^{-\theta,p}(\xi)$,
$\mathbf{e}_n^{\theta^\perp,p}(\xi)$
and 
$\mathbf{e}_n^{-\theta^\perp,p}(\xi)$.

\begin{figure}[t]
	\begin{picture}(10,160)(10,-80)
	\scalebox{.9}{
		\rput(1,0.5){\normalsize{\textcolor{blue}{$\theta$}}}
		\rput(-1,0.5){\normalsize{\textcolor{red}{$\theta^\perp$}}}
		\psline[linecolor=blue,linewidth=3pt]{->}(0,0)(1,1)
		\psline[linecolor=red,linewidth=3pt]{->}(0,0)(-1,1)
		\psline[linecolor=black,linewidth=1pt](-3,-3)(3,3)
		\psline[linecolor=black,linewidth=1pt](3,-3)(-3,3)
		\rput(0.2,-2.2){
			\scalebox{0.9}{
				\begin{minipage}{3cm}
				$\begin{psmatrix}[rowsep=0.4cm]
				\xi\in H_{\theta}\cap H_{\theta^\perp}
				\\[-3ex]
				\Downarrow
				\\[-2.25ex]
				\psframebox{\mathbf{e}_n^{\theta,p}(\xi),\mathbf{e}_n^{\theta^\perp,p}(\xi)}
				\\
				\psframebox{\widehat{\mathrm{E}^\theta}_n(\xi)e_i}
				\psset{nodesep=1pt}
				\ncline{->}{3,1}{4,1}{\footnotesize\trput{\raisebox{0.0cm}{\eqref{errorAK>1}}}}
				\end{psmatrix}
				$
				\end{minipage}
			}
		}
		\rput(0.,2.1){
			\scalebox{0.9}{
				\begin{minipage}{3cm}
				$\begin{psmatrix}[rowsep=0.4cm]
				\xi\in H_{-\theta}\cap H_{-\theta^\perp}
				\\[-3ex]
				\Downarrow
				\\[-2.25ex]
				\psframebox{\mathbf{e}_n^{-\theta,p}(\xi),\mathbf{e}_n^{-\theta^\perp,p}(\xi)}
				\\
				\psframebox{\widehat{\mathrm{E}^\theta}_n(\xi)e_i}
				\psset{nodesep=1pt}
				\ncline{->}{3,1}{4,1}{\footnotesize\trput{\raisebox{0.0cm}{\eqref{errorAK>1}}}}
				\end{psmatrix}
				$
				\end{minipage}
			}
		}
		\rput(2.75,0.2){
			\scalebox{0.9}{
				\begin{minipage}{3cm}
				$\begin{psmatrix}[rowsep=0.4cm]
				\xi\in H_{-\theta}\cap H_{\theta^\perp}
				\\[-3ex]
				\Downarrow
				\\[-2.25ex]
				\psframebox{\mathbf{e}_n^{-\theta,p}(\xi),\mathbf{e}_n^{\theta^\perp,p}(\xi)}
				\\
				\psframebox{\widehat{\mathrm{E}^\theta}_n(\xi)e_i}
				\psset{nodesep=1pt}
				\ncline{->}{3,1}{4,1}{\footnotesize\trput{\raisebox{0.0cm}{\eqref{errorAK>1}}}}
				\end{psmatrix}
				$
				\end{minipage}
			}
		}
		\rput(-2.75,0.2){
			\scalebox{0.9}{
			\begin{minipage}{3cm}
			$\begin{psmatrix}[rowsep=0.4cm]
			\xi\in H_{\theta}\cap H_{-\theta^\perp}
			\\[-3ex]
			\Downarrow
			\\[-2.25ex]
			\psframebox{\mathbf{e}_n^{\theta,p}(\xi),\mathbf{e}_n^{-\theta^\perp,p}(\xi)}
			\\
			\psframebox{\widehat{\mathrm{E}^\theta}_n(\xi)e_i}
			\psset{nodesep=1pt}
			\ncline{->}{3,1}{4,1}{\footnotesize\trput{\raisebox{0.0cm}{\eqref{errorAK>1}}}}
			\end{psmatrix}
			$
			\end{minipage}
			}
		}
	}
	\end{picture}
	\caption{\label{cuadrantes}Construction of the Fourier transform of the error term at the $n$-th iteration for fixed angle $\theta$ when $K>1$.}
\end{figure}
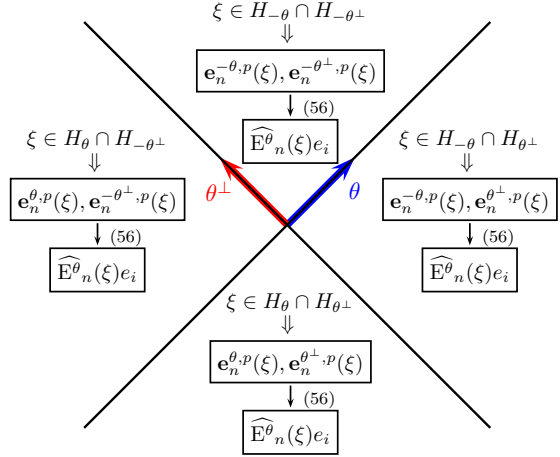
\begin{figure}[t]
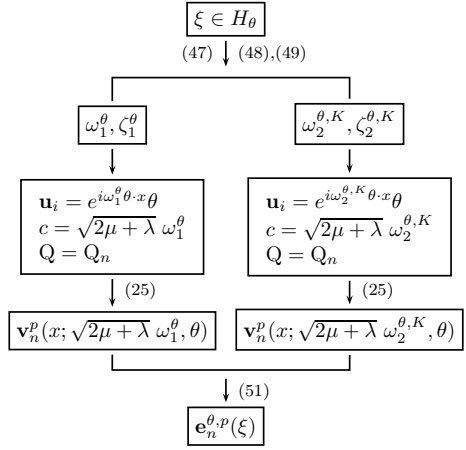

	\centerline{
		\scalebox{0.8}{
			\begin{minipage}{15cm}
				$$
				\begin{psmatrix}[colsep=-0.5cm,rowsep=0.5cm]
				&\psframebox{\xi\in H_\theta}&
				\\[-1ex]
				&{}&
				\\
				\psframebox{\omega_1^\theta,\zeta_1^\theta}&&\psframebox{\omega_2^{\theta,K},\zeta_2^{\theta,K}}
				\\
				\psframebox{
					\begin{array}{l}
					\mathbf{u}_i=e^{i\omega_1^\theta\theta\cdot x}\theta\\
					c=\sqrt{2\mu+\lambda}\ \omega_1^\theta\\
					\mathrm{Q}=\mathrm{Q}_n
					\end{array}}&&
				\psframebox{
					\begin{array}{l}
					\mathbf{u}_i=e^{i\omega_2^{\theta,K}\theta\cdot x}\theta\\
					c=\sqrt{2\mu+\lambda}\ \omega_2^{\theta,K}\\
					\mathrm{Q}=\mathrm{Q}_n
					\end{array}}
				\\
				\psframebox{\mathbf{v}^p_n(x;\sqrt{2\mu+\lambda}\ \omega_1^\theta,\theta)}&&
				\psframebox{\mathbf{v}^p_n(x;\sqrt{2\mu+\lambda}\ \omega_2^{\theta,K},\theta)}
				\\[-2ex]
				&{}&
				\\
				&\psframebox{\mathbf{e}^{\theta,p}_n(\xi)}&
				\psset{nodesep=2pt}
				\ncline{->}{1,2}{2,2}{\footnotesize\tlput{\raisebox{-0.2cm}{\eqref{w1z1}}}}{\footnotesize\trput{\raisebox{-0.2cm}{\eqref{w2},\eqref{z2}}}}
				\ncbar[angleA=90,nodesep=2pt,armA=0.3]{-}{3,1}{3,3}
				\ncline{->}{3,1}{4,1}
				\ncline{->}{3,3}{4,3}
				\ncline{->}{4,1}{5,1}{\footnotesize\trput{\raisebox{-0.2cm}{\eqref{LSeq}}}}
				\ncline{->}{4,3}{5,3}{\footnotesize\trput{\raisebox{-0.2cm}{\eqref{LSeq}}}}
				\ncbar[angleA=-90,nodesep=2pt,armA=0.25,armB=0.25cm
				]{-}{5,1}{5,3}
				\ncline{->}{5,2}{6,2}
				\ncline{->}{6,2}{7,2}{\footnotesize\trput{\raisebox{0.1cm}{\eqref{errorp}}}}
				\end{psmatrix}
				$$ 
			\end{minipage}
		}
	}
\caption{\label{diagrama_e_theta_p_n}Construction of the error $\mathbf{e}^{\theta,p}_n(\xi)$ at the n-th iteration when $\xi\in H_\theta$.}
\end{figure}
\section{Numerical approximation}
\label{Numerico}

In this section we describe the numerical method to find both the Born approximation and the sequence of loads $Q_n$ and $Q_n^\theta$ defined in the previous section from the corresponding scattering data. 
We focus on the backscattering situation given in Definition \ref{BornB} since the other
cases, namely the fixed angle scattering both for $K>1$ (Definition \ref{BornAK>1})  and $K<1$ (Definition \ref{BornAK<1}), can be treated similarly. 

We adapt the trigonometric collocation introduced in \cite{V} for the scalar equation, (Helmholtz equation). The required scattered waves are obtained using the numerical algorithm described in \cite{BC}. 

\subsection{Finite dimensional trigonometric space}

Given $R>0$, we define
$$
G_R=\left\{ x=(x_1,x_2)\in \R^2 \; : \; |x_k|<R,\; k=1,2 \right\} .
$$
The family of exponentials
$$
\varphi_j(x)=\frac{e^{i\pi j \cdot x/R}}{2R},  \quad j=(j_1,j_2)\in\Z^2,
$$
constitutes an orthonormal basis on $L^2(G_R)$ with the norm
$$\| u \|_0^2=\int_{G_R}|u(x)|^2 dx.$$ We also introduce the space $H^\eta=H^\eta(G_R)$ which consists of $2R-$multiperiodic functions (distributions) having finite the norm
$$
\| u \|^2_\eta=\sum_{j\in\Z^2} (1+|j|)^{2\eta} | \widehat{u}_j|^2  ,
$$
where
$$
\widehat{u}_j=\int_{G_R} u(x)\overline{\varphi_j(x)} dx=\frac{(u\chi_{G_R})\,^{\widehat{}} (\xi_j)}{2R}, \qquad \xi_j=\frac{j}{2R}, \quad j\in \Z^2.
$$
Note that $\widehat{u}_j$ is the j-th Fourier coefficient of $u$ with respect to the orthonormal family $\{\varphi_j\}_{j \in \mathbb{Z}^2}$.

We now introduce a finite dimensional approximation of $H^\eta$. Let us consider $h=2R/N$ with $N\in \N$ and a mesh on $G_R$ with grid points $jh$, $j\in \Z_h^2$ and
$$
\Z_h^2=\left\{ j=(j_1,j_2)\in \Z^2 \; : \; -\frac{N}{2} \leq j_k < \frac{N}{2}, \; k=1,2 \right\}.
$$
We denote by $\mathcal{T}_h$ the finite dimensional subspace of trigonometric polynomials of the form
$$
v^h=\sum_{j\in \Z_h^2} c_j \varphi_j, \qquad c_j\in \C.
$$
Any $v^h\in \mathcal{T}_h$ can be represented either through the Fourier coefficients
$$
v^h(x)=\sum_{j\in \Z_h^2}  \widehat{v^h}_j \; \varphi_j(x),
$$
or the nodal values
$$
v^h(x)=\sum_{j\in \mathbb{Z}_h^2} v^h(j h) \; \varphi^h_j(x),
$$
where $\varphi^h_j(kh)=\delta_{jk}$, more specifically
$$
\varphi^h_j(x)=\frac{h^2}{(2R)^2}\sum_{k\in \Z_h^2} e^{i\pi k \cdot (x-jh)/R}.
$$
For a given $v^h \in \mathcal{T}_h$, abusing notation, we write its Fourier coefficients as $\widehat{v^h}=\{\widehat{v^h}_j\}_{j\in\mathbb{Z}_h^2}$ and
its nodal values as
$v^h=\{v^h(jh)\}_{j\in\mathbb{Z}_h^2}$.
These Fourier coefficients and nodal values are related by the discrete Fourier transform $\mathcal{F}_h$ as follows:
\begin{equation*}\label{relation_F_F-1}
2 R\, \widehat{ v^h} =h^2\mathcal{F}_h v^h, \qquad v^h =\frac{2R}{h^2}\mathcal{F}_h^{-1} \widehat{ v^h}.
\end{equation*}

Here, as usual, $\mathcal{F}_h$ relates two sequences $x=\{x(n)\}_{n\in\mathbb{Z}_h^2}$ and $X=\{X(j)\}_{j\in\mathbb{Z}_h^2}$ according to
$$
\mathcal{F}_hx(j)=X(j)=\sum_{n\in\mathbb{Z}_h^2} x(n)e^{-i2\pi n/N}, \qquad j\in\mathbb{Z}_h^2, $$
and 
$$
\mathcal{F}_h^{-1}X(n)=x(n)= \frac{(2R)^2}{h^2}\sum_{j\in\mathbb{Z}_h^2} X(j)e^{i2\pi j/N}, \qquad n\in\mathbb{Z}_h^2.
$$
This definition coincides with the usual one in numerical codes (as MATLAB) up to a translation, since it considers $-\frac{N}{2}\le j_k<\frac{N}{2}$ instead of $0\le j_k<N$. This must be taken into account in the implementation.

The space $H^\eta$ and its finite dimensional approximation contains scalar functions but we need to extend these to vector-valued functions and to matrix-valued functions. This can be done component by component, so that
$$
\mathbf{v}^h=({v_1}^h,{v_2}^h) \in \mathcal{T}_h \times \mathcal{T}_h \quad\text{and}\quad  
\mathrm{Q}^h    =\begin{pmatrix}
{Q_{11}}^h & {Q_{12}}^h \\ {Q_{21}}^h & {Q_{22}}^h
\end{pmatrix} \in \mathcal{M}_{2\times 2}(\mathcal{T}_h).
$$ 



\subsection{Finite dimensional setting}

We first focus on the numerical method to approximate $\mathrm{Q}_B$ from the backscattering data, that we assume known.

Using the definition given in (\ref{Def_BornB}) we can construct the following finite dimensional version of $\mathrm{Q}_B$: 
\begin{equation}
    \label{Def_BornB_h_0}
    \mathrm{Q}_{B}^h(x)= \sum_{j \in \mathbb{Z}_h^2}\mathrm{Q}_{B}^h(jh)\varphi^h_j(x) \in \mathcal{M}_{2\times 2}(\mathcal{T}_h),
\end{equation}
where for $i=1,2$
\begin{equation}
    \label{Def_BornB_h}
    h^2 \mathrm{Q}_B^he_{i}=2R \mathcal{F}_h^{-1}  \widehat{\mathrm{Q}^h_B} e_i,
\end{equation}
with
$$
\widehat{\mathrm{Q}^h_B} e_i=\frac{1}{2R}
\left\{   (e_{i} \cdot \theta(\xi_j))\mathbf{v}_{\infty}^{p}(\omega(\xi_j),\theta(\xi_j))     
 +
( e_{i}\cdot \theta^\perp(\xi_j))\mathbf{v}_{\infty}^{s}(\omega(\xi_j),\theta(\xi_j))  \right\}_{j \in \mathbb{Z}^2_h  },
$$
$\xi_j=\frac{j}{2R}$, and $\omega(\xi_j)$, $\theta(\xi_j)$ given by 
\begin{equation*}
\omega(\xi_j)=\frac{|\xi_j|}{2}, \quad \theta(\xi_j)=-\frac{\xi_j}{|\xi_j|},
\end{equation*}
in such a way that $\xi_j=-2\omega(\xi_j)\theta(\xi_j) $ (see Remark \ref{OmegaTheta} above). 

Note that $\mathrm{Q}_B^h e_i$ is computed from a single inverse discrete Fourier transform from the values of the far fields $\mathbf{v}_\infty^p(\omega(\cdot),\theta(\cdot))$ and $\mathbf{v}_\infty^s(\omega(\cdot),\theta(\cdot))$ at the mesh points $\xi_j$ with $j\in \Z_h^2$. More precisely, the numerical approximation is obtained from the following process:

\bigskip

{\bf Algorithm 2:}
\begin{enumerate}
\item  Choose $h$ according to the mesh grid where we will compute the nodal values of $\mathrm{Q}_B^h$.
\item Construct the mesh $\xi_j=j/(2R)$ with $j\in\Z_h^2$.
\item Evaluate $\mathbf{v}_\infty^p(\omega(\cdot),\theta(\cdot))$ and $\mathbf{v}_\infty^s(\omega(\cdot),\theta(\cdot))$ given respectively by (\ref{b_p}) and (\ref{b_s}) at the mesh points $\xi_j$.
\item Invert the discrete Fourier transform in (\ref{Def_BornB_h}) to obtain the values of both $\mathrm{Q}_B^h e_1$ and $\mathrm{Q}_B^h e_2$ at the nodes $x_j$. This gives the four componentes of the matrix $\mathrm{Q}_B^h$.
\end{enumerate}

We now consider the iterative process described in \eqref{algoritmo_back}. The following algorithm is a detailed discrete version of the Algorithm 1 in the introduction.
\bigskip

{\bf Algorithm 3:}
\begin{enumerate}
\item  Choose $h$ according to the mesh grid where we will compute the nodal values of $\mathrm{Q}_B^h$.
\item Construct the mesh $\xi_j=j/(2R)$ with $j\in\Z_h^2$.
\item Compute $\mathrm{Q}_B^h$ following the algorithm 2 above and define $\mathrm{Q}_1=\chi_{R}\mathrm{Q}_B^h$.
\item We set $M \in \mathbb{N}$  and for $n=1,2, \cdot \cdot \cdot ,M$ we define $\mathrm{Q}_{n+1}$ from $\mathrm{Q}_n$ as follows:
\begin{enumerate}
    \item[(i)] 
    Solve (\ref{LSeq}) four times with $\mathrm{Q}=\mathrm{Q}_n$ always but with four different incident waves $\mathrm{u}_i$ to get the four corresponding scattered solutions (see Figure \ref{diagramaEn} for more details).
    To do this use the method given in \cite{BC}.
    \item [(ii)] Evaluate $\mathbf{e}^p(\omega(\cdot),\theta(\cdot))$ and $\mathbf{e}^s(\omega(\cdot),\theta(\cdot))$ given respectively by (\ref{hp}) and (\ref{hs}) at the mesh points $\xi_j$.
    \item [(iii)] Construct the finite dimensional version of $\mathrm{E}_n$
    $$
    \mathrm{E}_{n}^h(x)= \sum_{j \in \mathbb{Z}_h^2}\mathrm{E}_{n}^h(j h)\varphi^h_j(x) \in \mathcal{M}_{2\times 2}(\mathcal{T}_h)
    $$
    where for $i=1,2$
$$ h^2 \mathrm{E}_n^he_{i}=2R \mathcal{F}_h^{-1}  \widehat{E_n^h}e_i $$
with
$$ \widehat{E_n^h}e_i=\frac{1}{2R}\left\{   (e_{i} \cdot \theta(\xi_j))\mathbf{e}^{p}(\omega(\xi_j),\theta(\xi_j))     
 +
( e_{i}\cdot \theta^\perp(\xi_j))\mathbf{e}^{s}(\omega(\xi_j),\theta(\xi_j))  \right\}_{j \in \mathbb{Z}^2_h  }$$
    \item[(iv)]  Write
$$\mathrm{Q}_{n+1}^h=\chi_{R}\mathrm{Q}_B^h-\chi_{R}\mathrm{E}_n^h.$$

\end{enumerate}
\end{enumerate}

\subsection{Convergence of the numerical approximation}

Here we give estimates for the error of the
numerical approximation ${\mathrm{Q}_B}_h$ in (\ref{Def_BornB_h_0}) with respect to
a periodized version of the Born approximation $\mathrm{Q}_{B}$. The following result is an easy consequence of the scalar analogous in \cite{BCR} and it establishes the convergence of the numerical Born approximation to the continuous one, up to a possible aliasing effect. 

\begin{theorem} \label{th_1}
Let $\mathrm{Q}_B$ the Born approximation of a load $\mathrm{Q}(x)$ defined by
(\ref{Def_BornB}). Let $\mathrm{Q}_B^\sharp$ the periodized version of $\mathrm{Q}_B$ defined as
\begin{equation}\label{eq_52e}
 \mathrm{Q}_B^\sharp(x)=\sum_{j\in \Z^2_h} \mathrm{Q}_B(x+2R j).
\end{equation}
If $\mathrm{Q}_B^\sharp \in H^\eta$ for some $\eta>0$, then
\begin{equation*} 
\| \mathrm{Q}_B^h -\mathrm{Q}_B^\sharp \|_0 \leq h^\eta \| \mathrm{Q}_B^\sharp
\|_\eta, \quad h>0,
\end{equation*}
where $\mathrm{Q}_B^h\in \mathcal{M}_{2\times 2}(\mathcal{T}_h)$ is the solution of (\ref{Def_BornB_h_0}) and $\|\mathrm{Q}\|_\eta = \max_{i,j=1,2}\|{\mathrm{Q}_{ij}} \|_\eta$.
\end{theorem}

\begin{remark}
The periodic version of $\mathrm{Q}_B$ given in (\ref{eq_52e}) will
coincide with $\mathrm{Q}_B$ in $G_R$ only if $\mathrm{Q}_B$ is compactly supported
in $G_R$. We do not know if this is the case, in general, but the numerical experiments below suggest that this is not the case. 
\end{remark}

\section{Numerical experiments}
\label{experimentos}

In this section we illustrate with numerical experiments the Born approximation both for backscattering and fixed angle data. We also show the efficiency of the iterative algorithms described  in Section \ref{SeccionAlgoritmos}  that we discretize following Section 5.  

\subsection{Born approximation}
We focus on the backscattering case since we did not find significant differences when simulating the Born approximation for fixed angle data.

We consider matrices of the form 
\begin{equation} \label{eq:matpot}
\mathrm{Q}(x) = q(x) \left( 
\begin{array}{ll} 1&1 \\ 1&1 
\end{array}
\right), \qquad x=(x_1,x_2),
\end{equation}
for a scalar function $q(x)$. We have tried with other situations as diagonal, non-symmetric and anti-symmetric matrices with similar results.  We present here two different situations corresponding to a discontinuous load given by 
\begin{equation} \label{eq:pot1}
q(x)=\left\{ 
\begin{array}{ll}
1.2, & \mbox{ if $|x_1|+|x_2|<0.2$}\\
1, & \mbox{ if $0.6< |x|<0.8$}\\
0, & \text{otherwise}.
\end{array}
\right.
\end{equation}
and a smooth one given by 
\begin{eqnarray} \nonumber
	q(x)&=&\max(0,e^{-5|x-(0.5,0)|^2})+1.5e^{-4|x-(-1,0.8)/2|^2}\\&+&2e^{-7|x-0.4(-1,-1)|^2-0.4}. \label{eq:pot2}
\end{eqnarray}
Both functions have support in the domain $|x|<1$. However, as computational domain we choose the larger region $|x|< 2=R$ since this allows us to recover the load from the inverse Fourier transform in a finer mesh for the frequency space, and this improves the precision of the reconstruction.    

In the experiments below we take the Lamé parameters $\lambda=2$ and $\mu=1$ but other values of the Lamé parameters provide similar results. 

It is worth mentioning that, as we said in the introduction, the known theoretical results for the existence of the Born approximation require the matrix $\mathrm{Q}$ to be symmetric and with all the components $C^1$ functions. However, as we show below none of these conditions seem to be relevant when simulating the Born approximation, and the known theoretical results can be probably extended to more general situations.

We also note that the Born approximation is defined as a complex function even if the original matrix is real. In our experiments we have observed that the Born approximation has a small complex component that we do not consider. More precisely, we compare the matrix load with the real part of its Born approximation.    

\bigskip

\begin{figure}
    \centering
    \psfrag{x}{\footnotesize{$x_1$}}
    \psfrag{poten}{\footnotesize{$\mathrm{Q}_{11}$}}
    \psfrag{approx}{\footnotesize{${\mathrm{Q}_B}_{11}$}}
    \begin{tabular}{cc}
    \includegraphics[width=6cm]{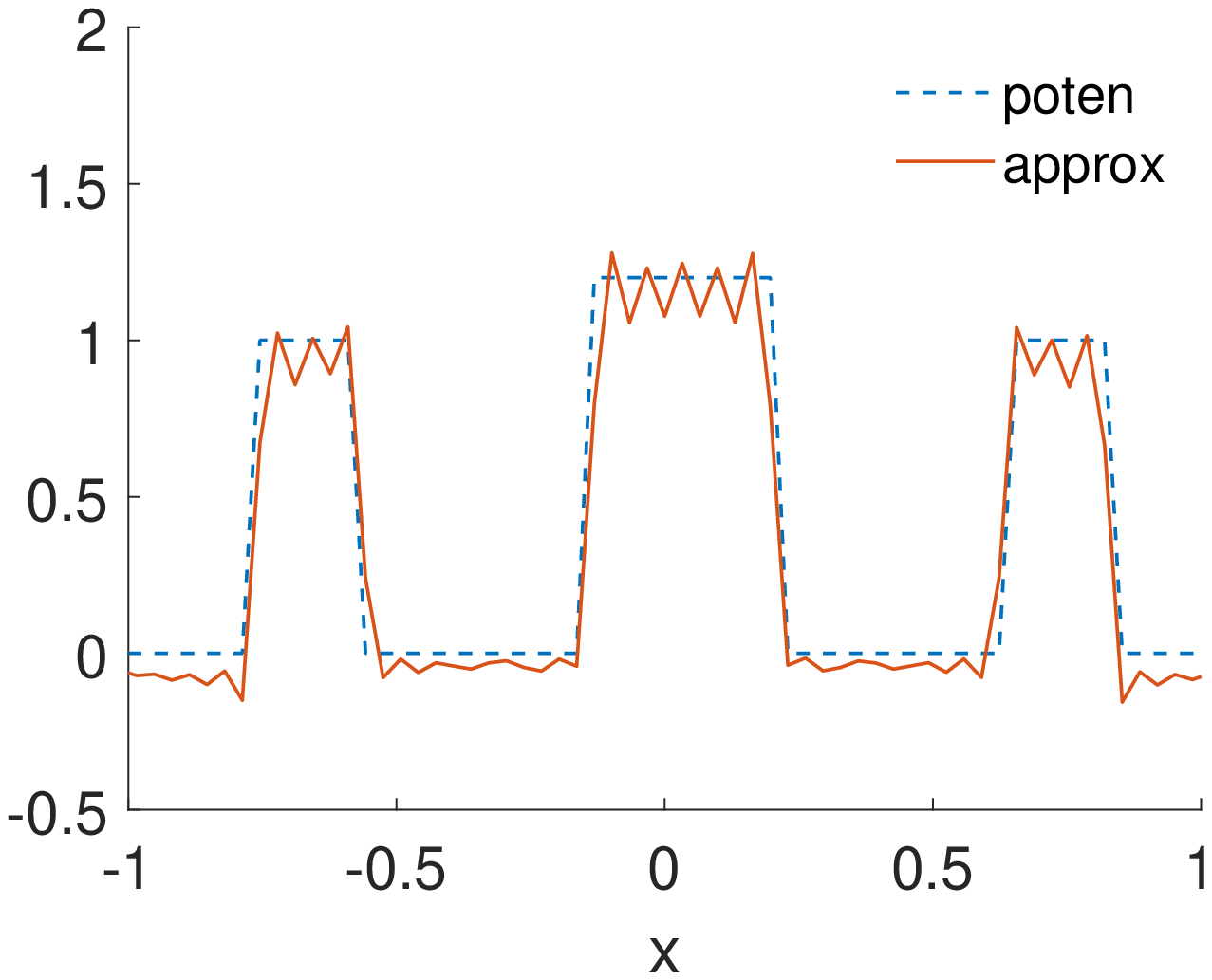}
    &
    \includegraphics[width=6cm]{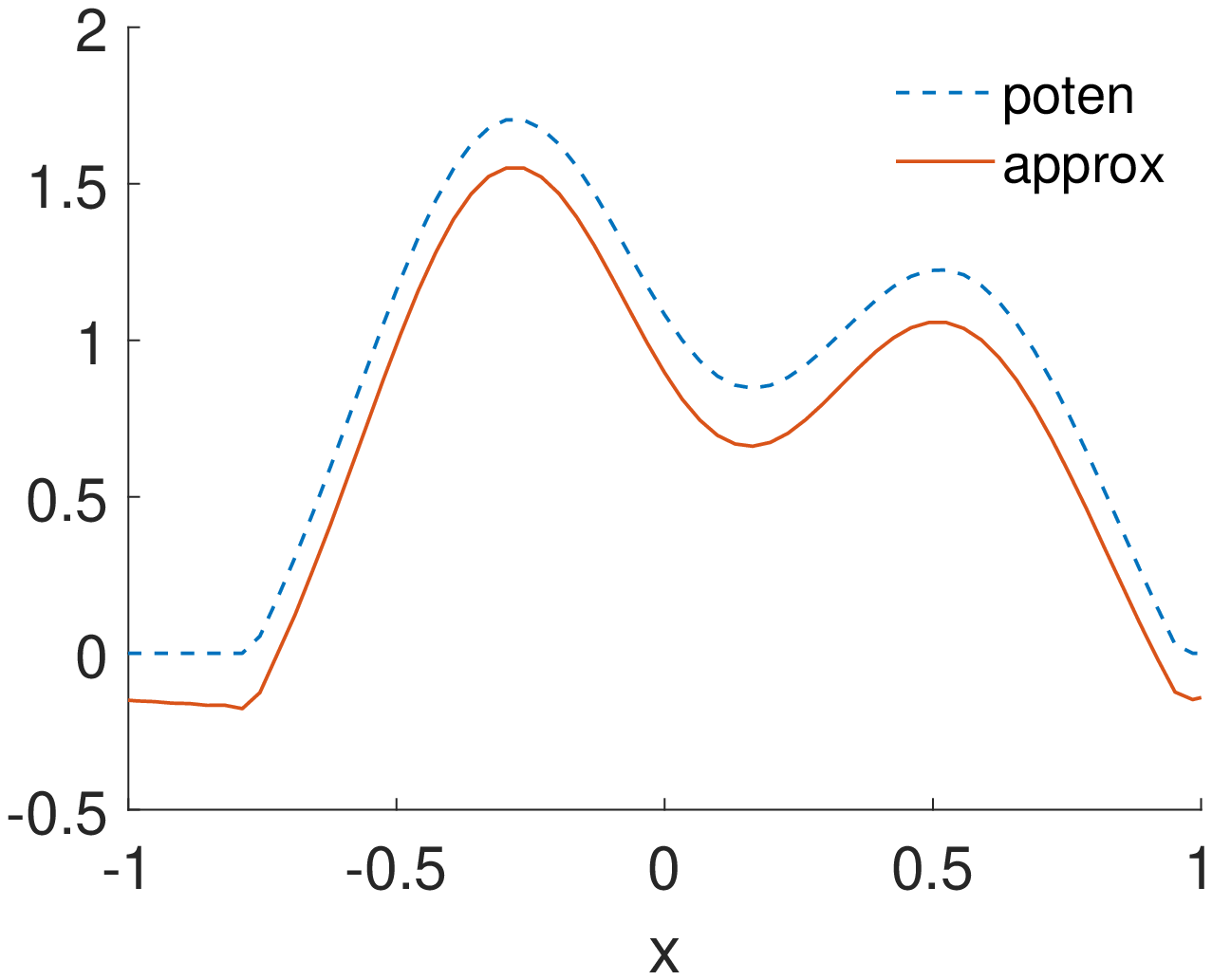} 
    \end{tabular}
    \caption{Experiment 1: central section of the first component of the load (dash) and the real part of its Born approximation (solid) with backscattering data. The discretization parameter is $h=2^{-6}$.}
    \label{fig:ex1}
\end{figure}

{\bf Experiment 1.} We first illustrate how good is the Born approximation to simulate the matrix load. In Figure \ref{fig:ex1} we compare the central section of both the first component of the load $\mathrm{Q}_{11}(x_1,0)$ and the real part of its Born approximation ${\mathrm{Q}_B}_{11}(x_1,0)$. We observe that the Born approximation provides a fairly good approximation of the load, even if it is not smooth. 

\bigskip

\begin{figure}
    \centering
    \psfrag{x}{\footnotesize{$x_1$}}
    \psfrag{poten}{\footnotesize{$\mathrm{Q}_{11}$}}
    \psfrag{approx}{\footnotesize{${\mathrm{Q}_B}_{11}$}}
    \begin{tabular}{cc}
    \includegraphics[width=6cm]{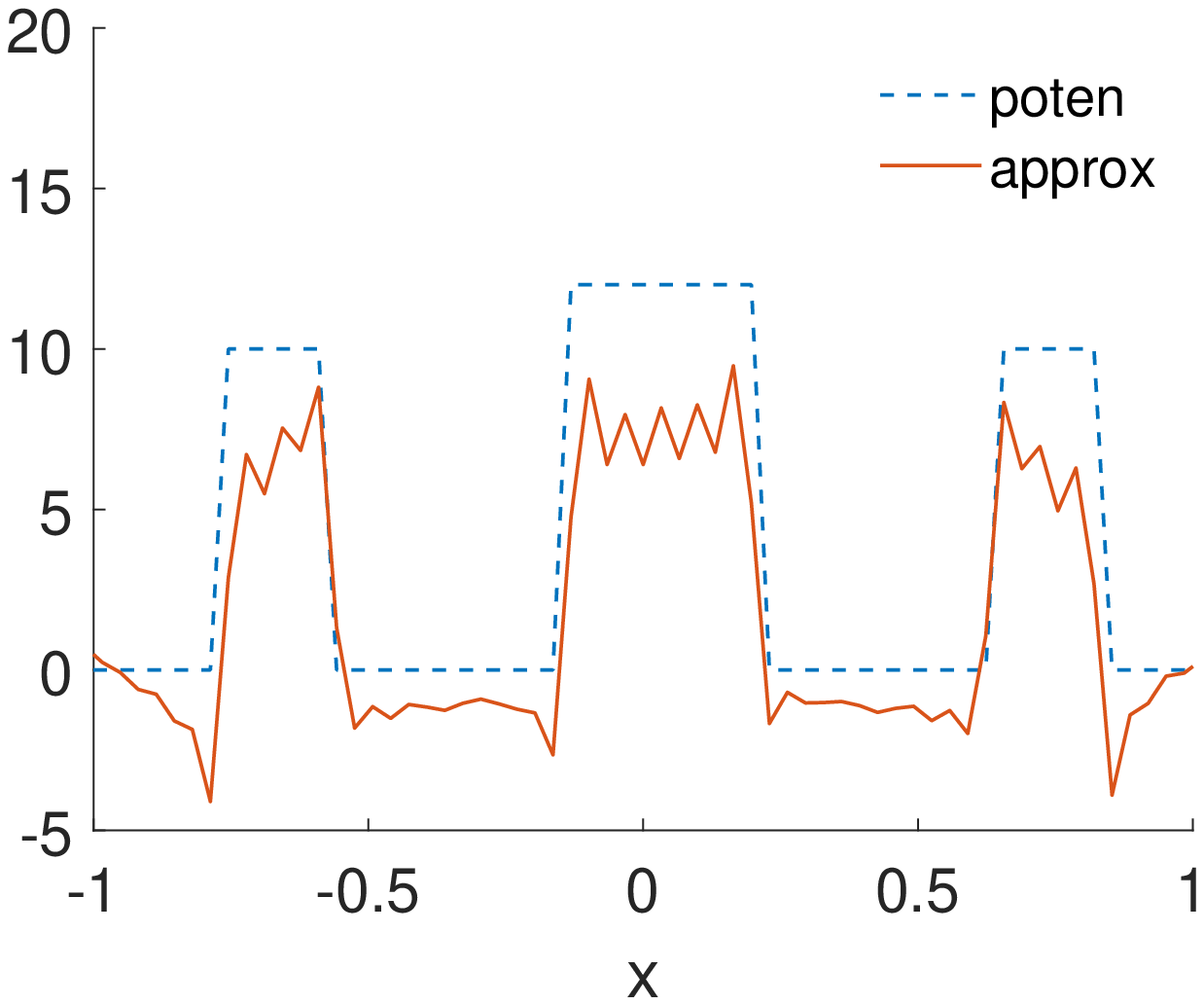} &
    \includegraphics[width=6cm]{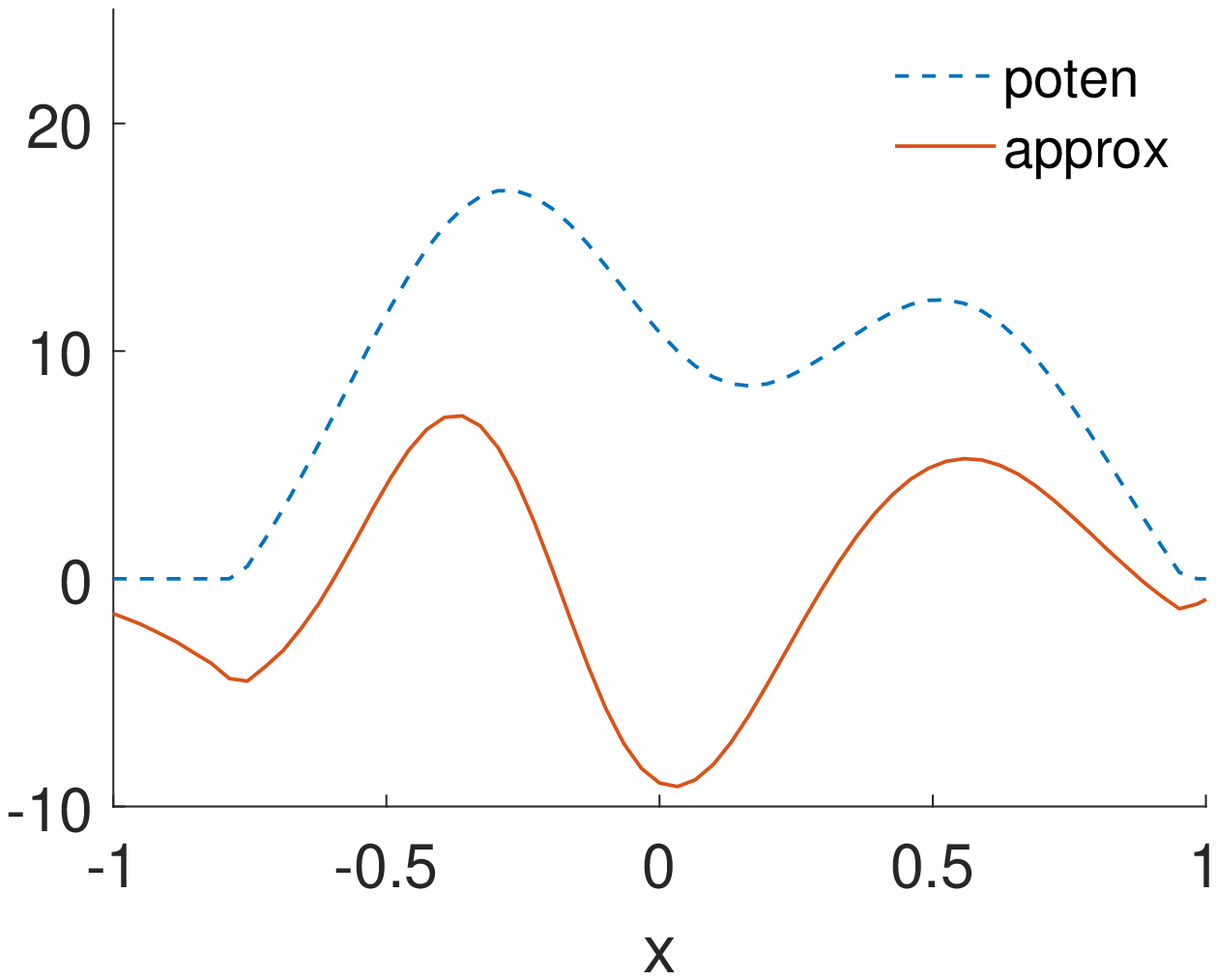} 
    \end{tabular}
    \caption{Experiment 2: central section of the first component of the load (dash) and the real part of its Born approximation (solid) obtained with backscattering data. The discretization parameter here is $h=2^{-6}$}
    \label{fig:ex2}
\end{figure}

{\bf Experiment 2.} Here we compute the Born approximation for larger loads. More precisely, we have multiplied by 10 both loads considered in the previous experiment. We see in Figure \ref{fig:ex2} that the Born approximation is not so close but still recover the aspect of the load. We also appreciate that the Born approximation is able to to detect the discontinuities in the left simulation. The recovery of singularities is a known property of the Born approximation both for backscattering and fixed angle scattering data, even for large loads (see \cite[Corollary 1.2]{bfprv1} and \cite[Theorem 1.1]{bfprv3}). We appreciate this better in the next experiment.  

\bigskip

\begin{figure}
    \centering
    \psfrag{x}{\footnotesize{$x_1$}}
    \psfrag{potenc}{\footnotesize{$\mathrm{Q}_{11}$}}
    \psfrag{Born}{\footnotesize{${\mathrm{Q}_B}_{11}$}}
    \begin{tabular}{cc}
    \includegraphics[width=5.5cm]{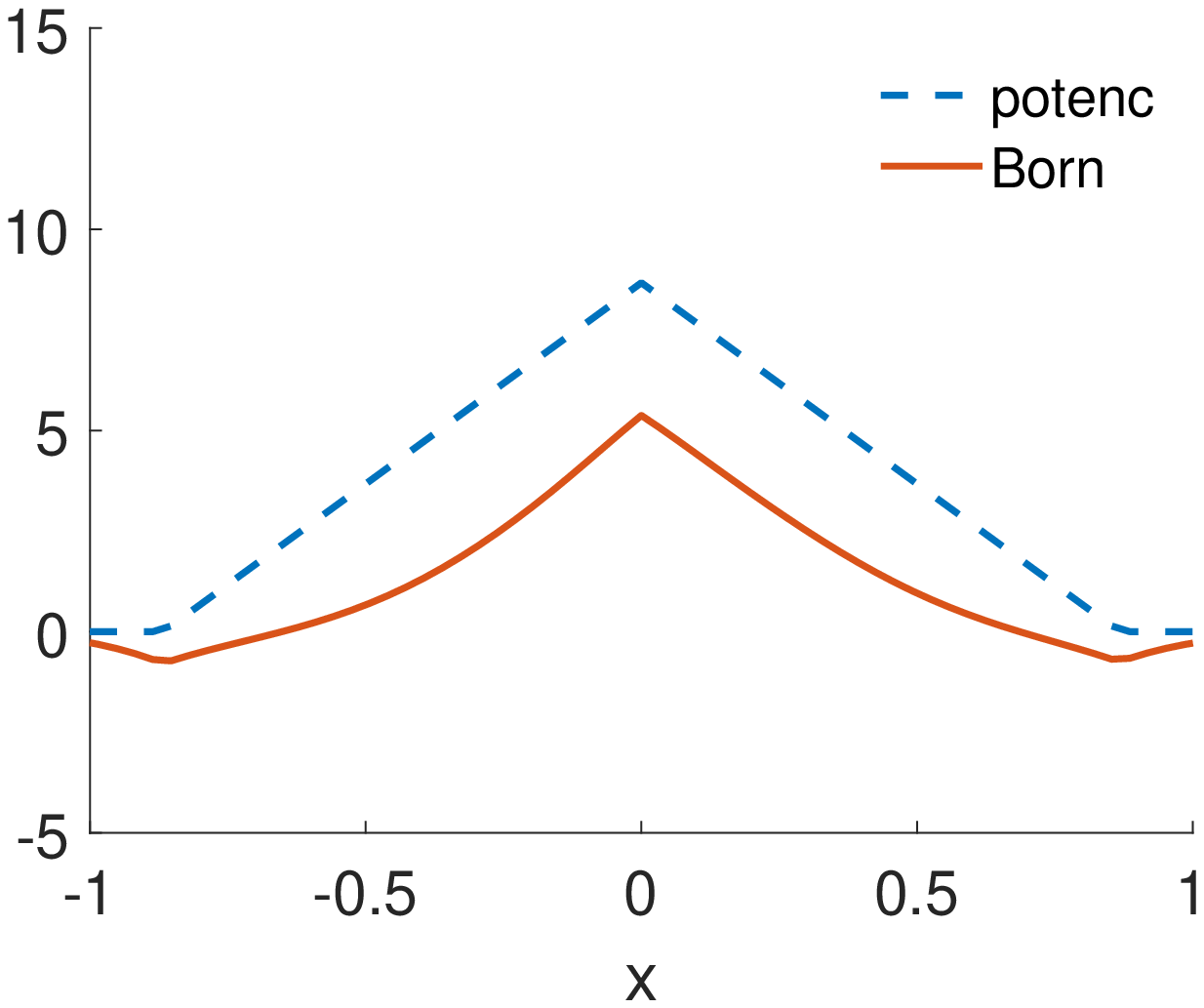} &
    \includegraphics[width=5.5cm]{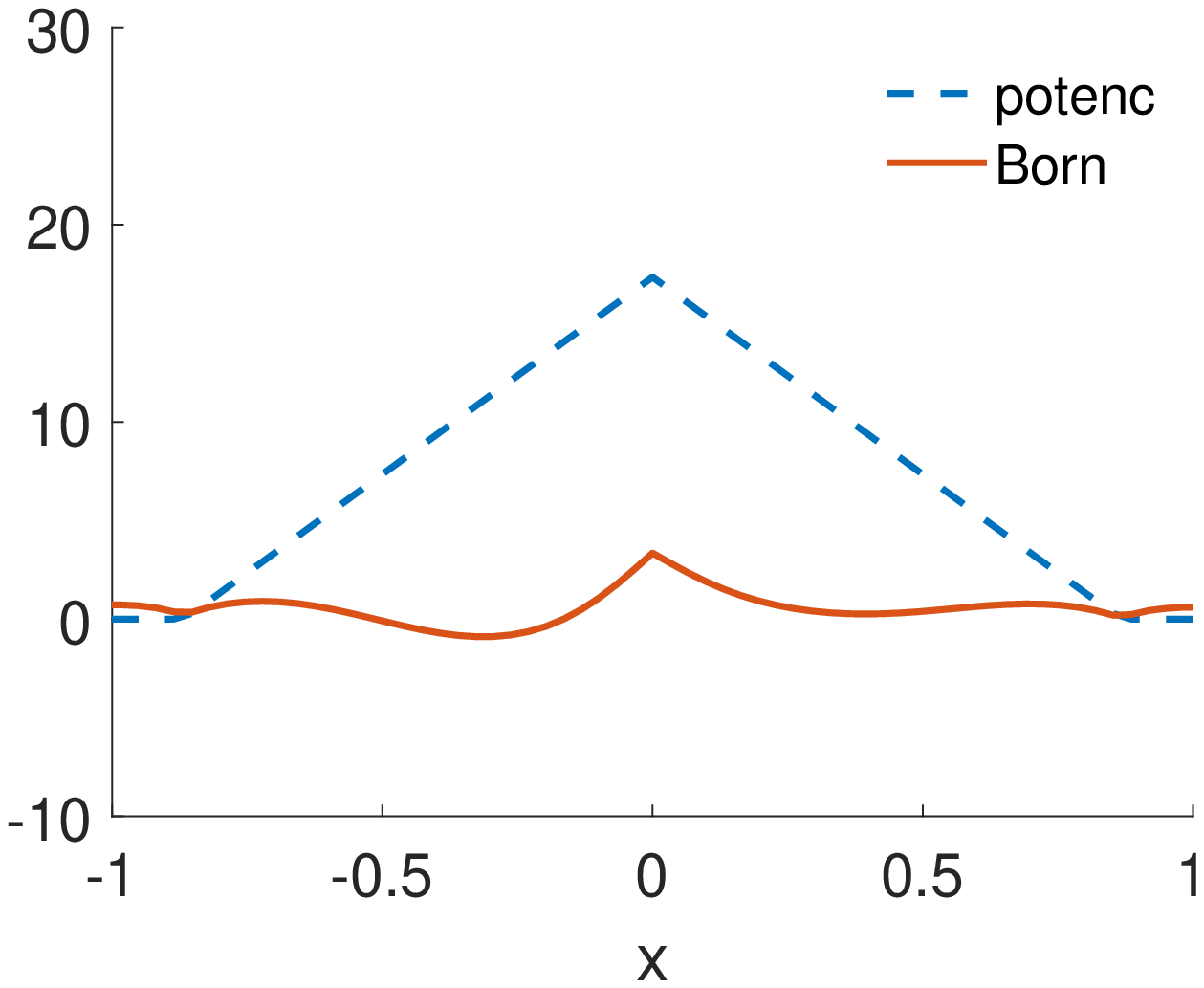} 
    \end{tabular}
    \caption{Experiment 3: central section of the first component of a Lipschitz load (dash) and its Born approximation (solid). The size of the load is twice larger in the right figure. We observe that the Born approximation may be a bad approximations when the load is large but still recovers the position of the singularities.} 
    \label{fig:ex4}
\end{figure}

{\bf Experiment 3.} Now we illustrate further the recovery of singularities of the Born approximation with a new experiment. We have considered the diagonal Lipschitz load given by $\mathrm{Q}(x)=\alpha (1-(|x_1|+|x_2|)) \mathrm{I}$ for different values of $\alpha$. In Figure \ref{fig:ex4} the central cross section of the load $\mathrm{Q}_{11}(x_1,0)$ and its Born approximation ${\mathrm{Q}_B}_{11}(x_1,0)$ for $\alpha=10$ (left) and $\alpha=20$ (right) are plotted. We observe that, even if the Born approximation is far from the load, it has a jump in the derivative at the same point as $\mathrm{Q}_{11}$ has it. 

\bigskip

\begin{figure}
    \centering
    \psfrag{x}{\footnotesize{$x_1$}}
    \begin{tabular}{cc}
    \includegraphics[width=6cm]{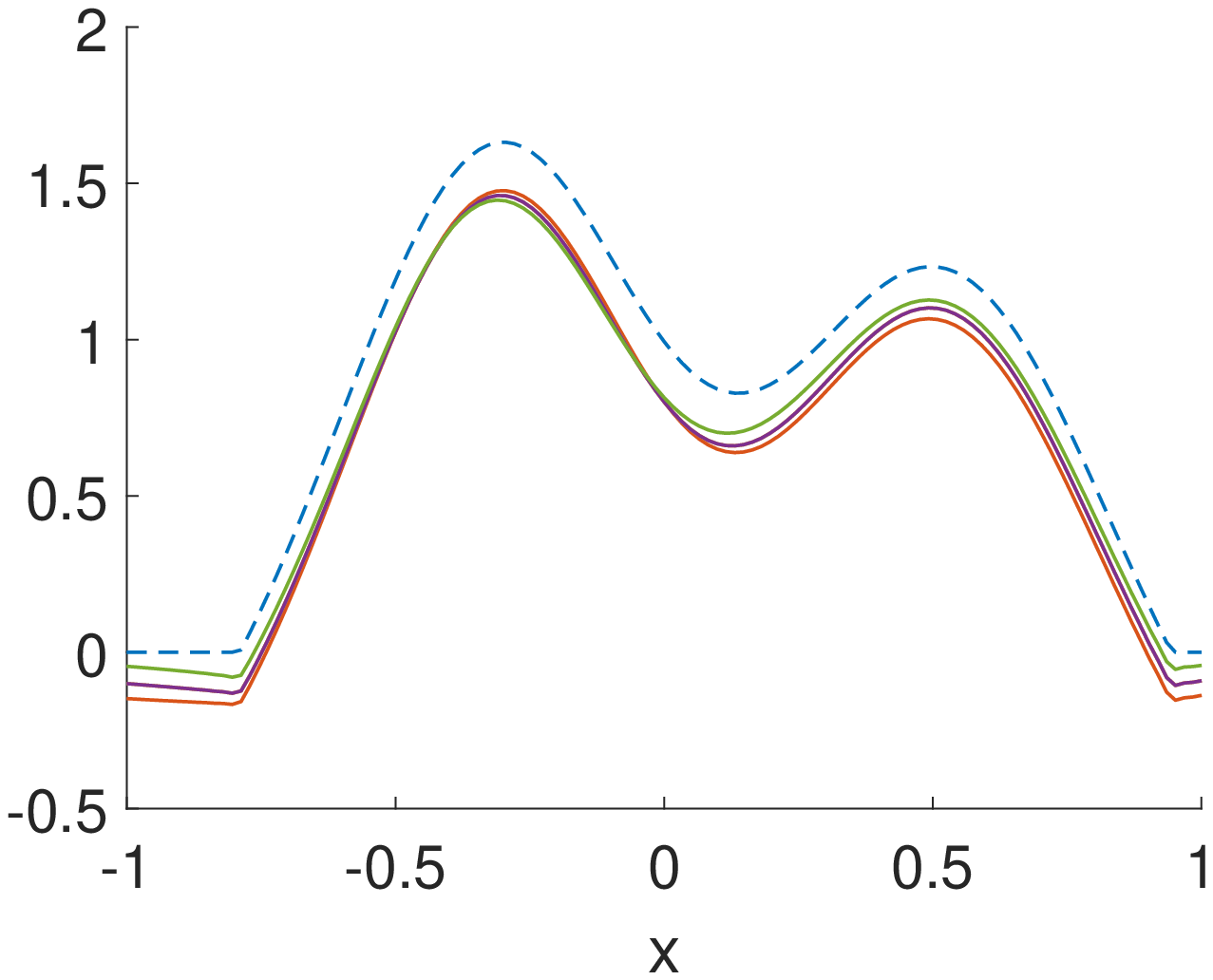} &
    \includegraphics[width=6cm]{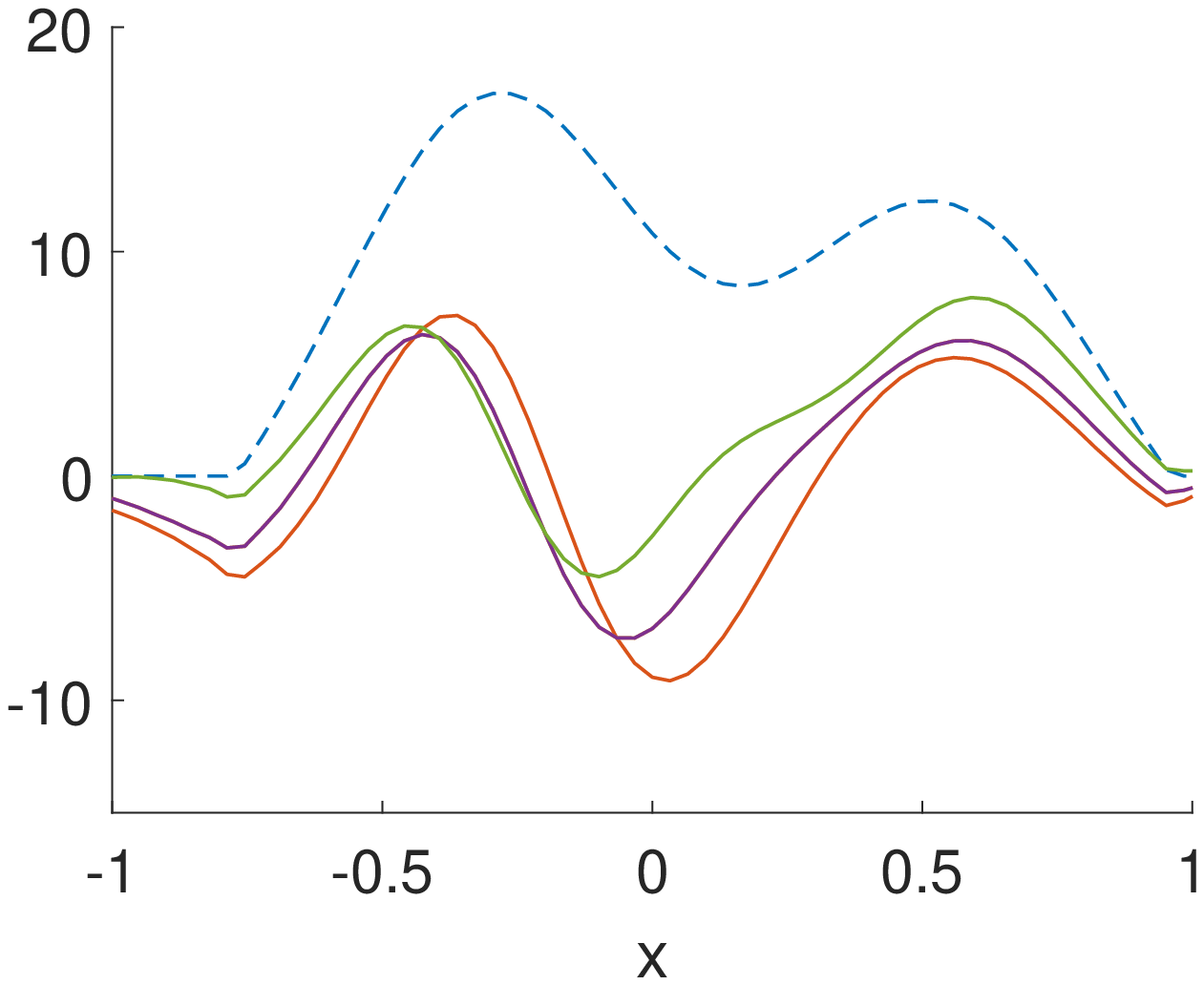} 
    \end{tabular}
    \caption{Experiment 3: Different components of the Born approximations ${\mathrm{Q}_B}_{\mathrm{ij}},\; i,j=1,2$ corresponding to a smooth load with equal components and $h=2^{-6}$. The difference is clearer in the right hand simulation which corresponds to a load 10 times larger than the left one.} 
    \label{fig:ex3}
\end{figure}

{\bf Experiment 4.}  In the experiment 1 above we assumed that the matrix load $\mathrm{Q}$ has the four components equal. However this is not the case for the Born approximation. In Figure \ref{fig:ex3} we plot the four components of the Born approximation matrix when the load is given by \eqref{eq:matpot}, with the smooth scalar function \eqref{eq:pot2}. We see that they are different, even if it is not the case in the original load. In fact, we observed that the Born approximation is a symmetric matrix, i.e. ${\mathrm{Q}_B}_{12}={\mathrm{Q}_B}_{21}$ but this function is different from ${\mathrm{Q}_B}_{11}$ and ${\mathrm{Q}_B}_{22}$.  

\bigskip



\subsection{Iterative algorithms}

We now check the performance of the iterative algorithms proposed in Section \ref{SeccionAlgoritmos} to approximate the load. We have added $5\%$ relative noise in the scattering data to simulate more realistic situations where measurements contain some errors. We observed that the behavior of the algorithms is quite robust to such noisy data.  

\bigskip

\begin{figure}
    \centering
    \psfrag{n}{\footnotesize{$n$}}
    \begin{tabular}{cc}
    \includegraphics[width=6cm]{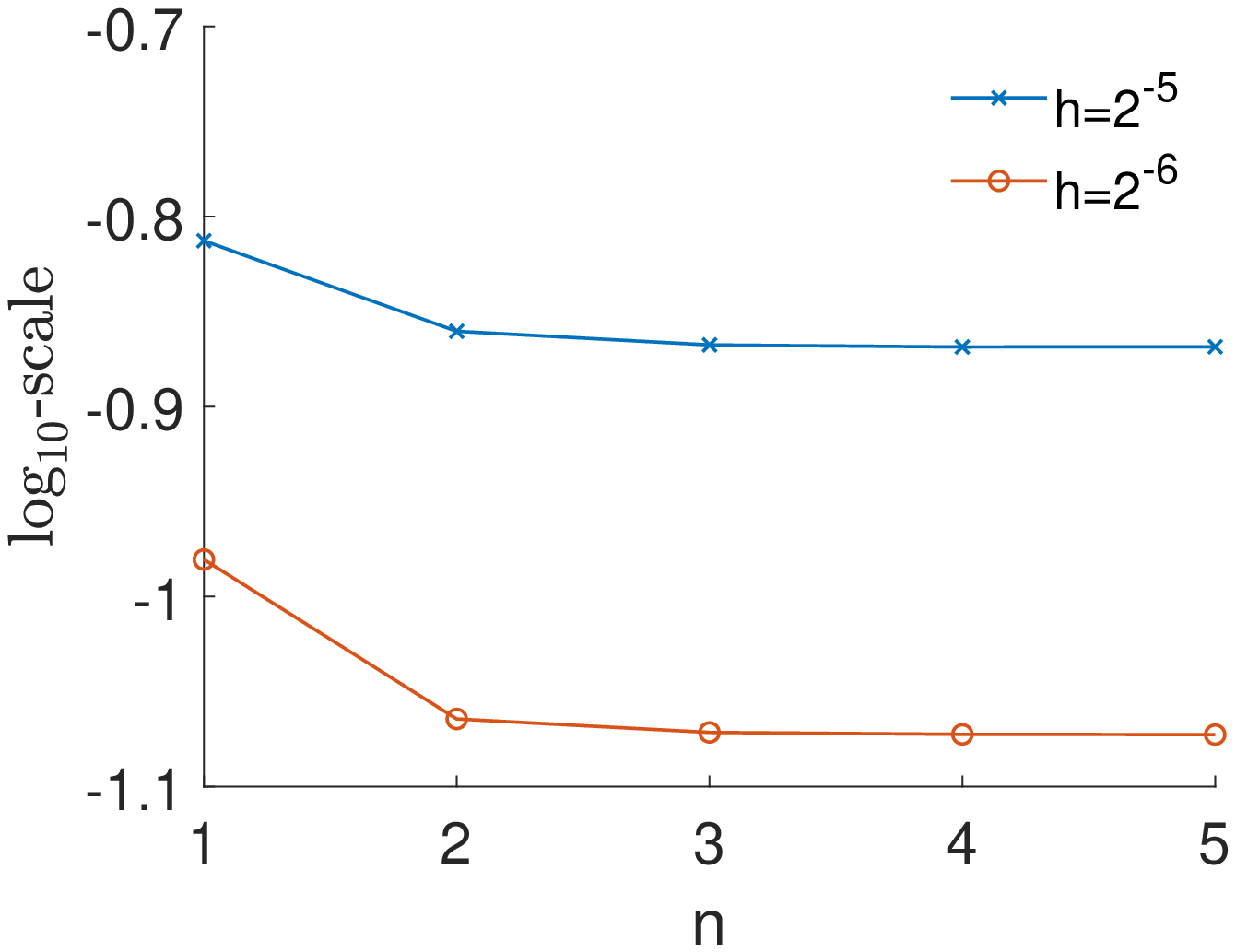} &
    \includegraphics[width=6cm]{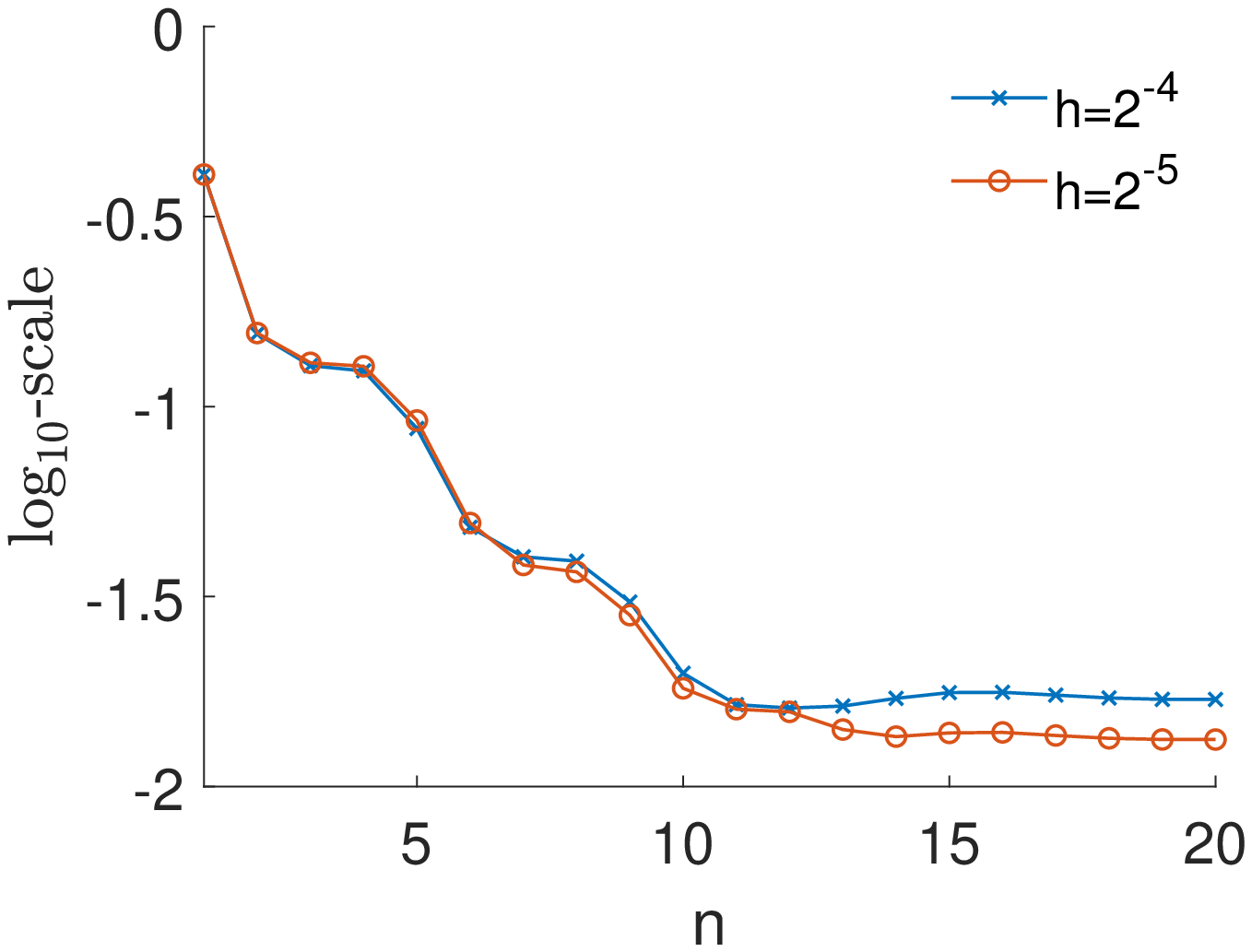} 
    \end{tabular}
    \caption{Experiment 5: error behavior in terms of the iterations for discontinuous (left) and smooth (right) loads.} 
    \label{fig:ex5}
\end{figure}

{\bf Experiment 5:} We first consider the backscatering case with the two loads given in \eqref{eq:pot1}-\eqref{eq:pot2}. We define the error at each iteration $n$ as
$$
\mbox{error}(n) = \max_{i,j=1,2} \left( h^2 \sum_{k,l\in \mathbb{Z}^{2}_h} \left| \mathrm{Q}_{ij}(x_{k},x_l)- Re \; (\mathrm{Q}_n)_{ij}
(x_{k},x_l) \right|^2 \right)^{1/2},
$$
where $Re$ denotes the real part. In Figure \ref{fig:ex5} we show the error behavior in terms of the number of iterations, both when $h=2^{-5}$ and $h=2^{-6}$. We observe that the discontinuous load is more sensitive to the discretization parameter $h$. This is natural due to the fact that we recover the matrix from its Fourier transform and this necessarily introduces a high frequency filtering that we appreciate in the oscillations appearing in the approximations (see also the left hand simulation in Figure \ref{fig:ex1}). A finer mesh provides a better approximation of more frequencies in the Fourier transform of the Born approximation. 

\begin{figure}
    \centering
    \psfrag{n}{\footnotesize{$n$}}
    \begin{tabular}{cc}
    \includegraphics[width=6cm]{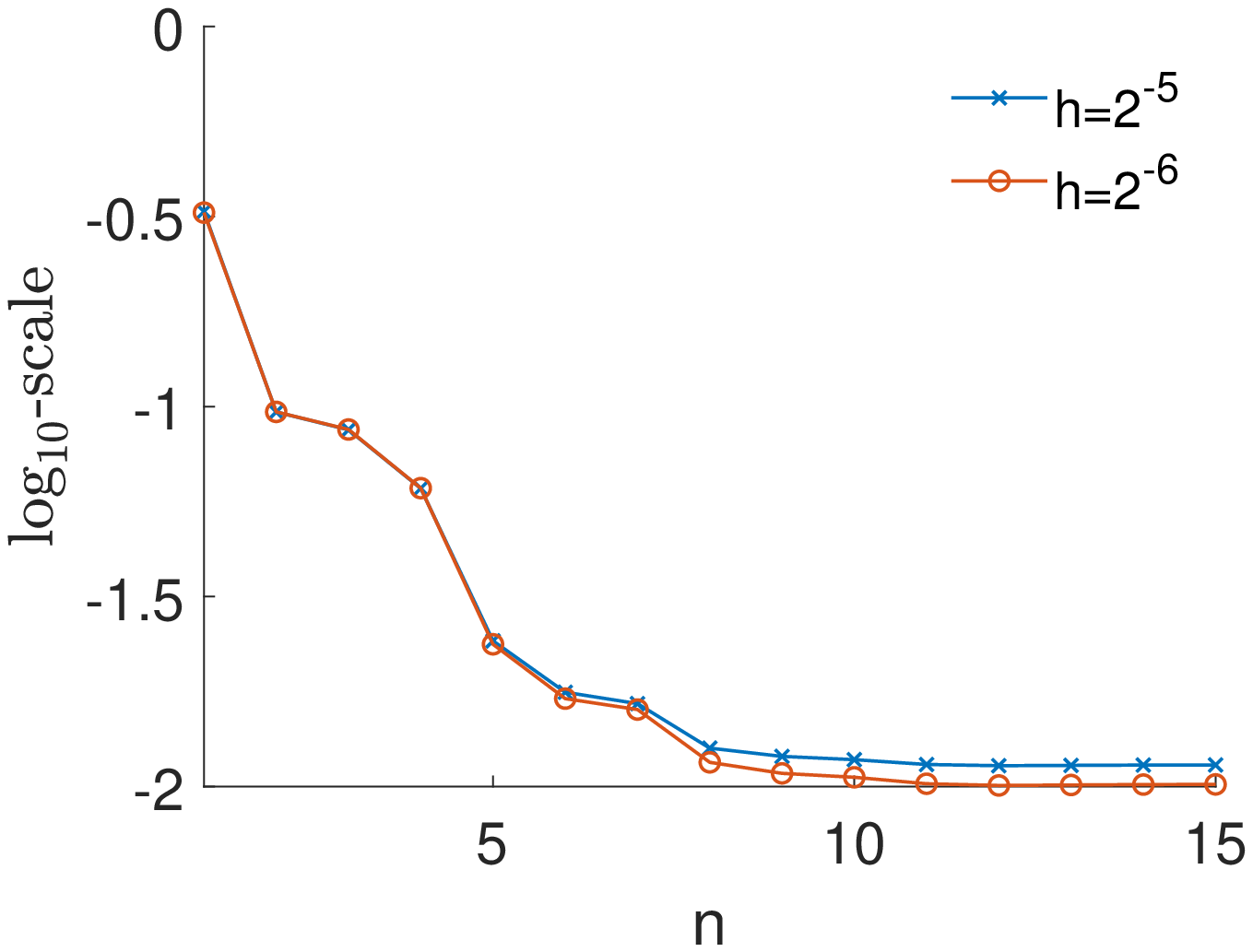} &
    \includegraphics[width=6cm]{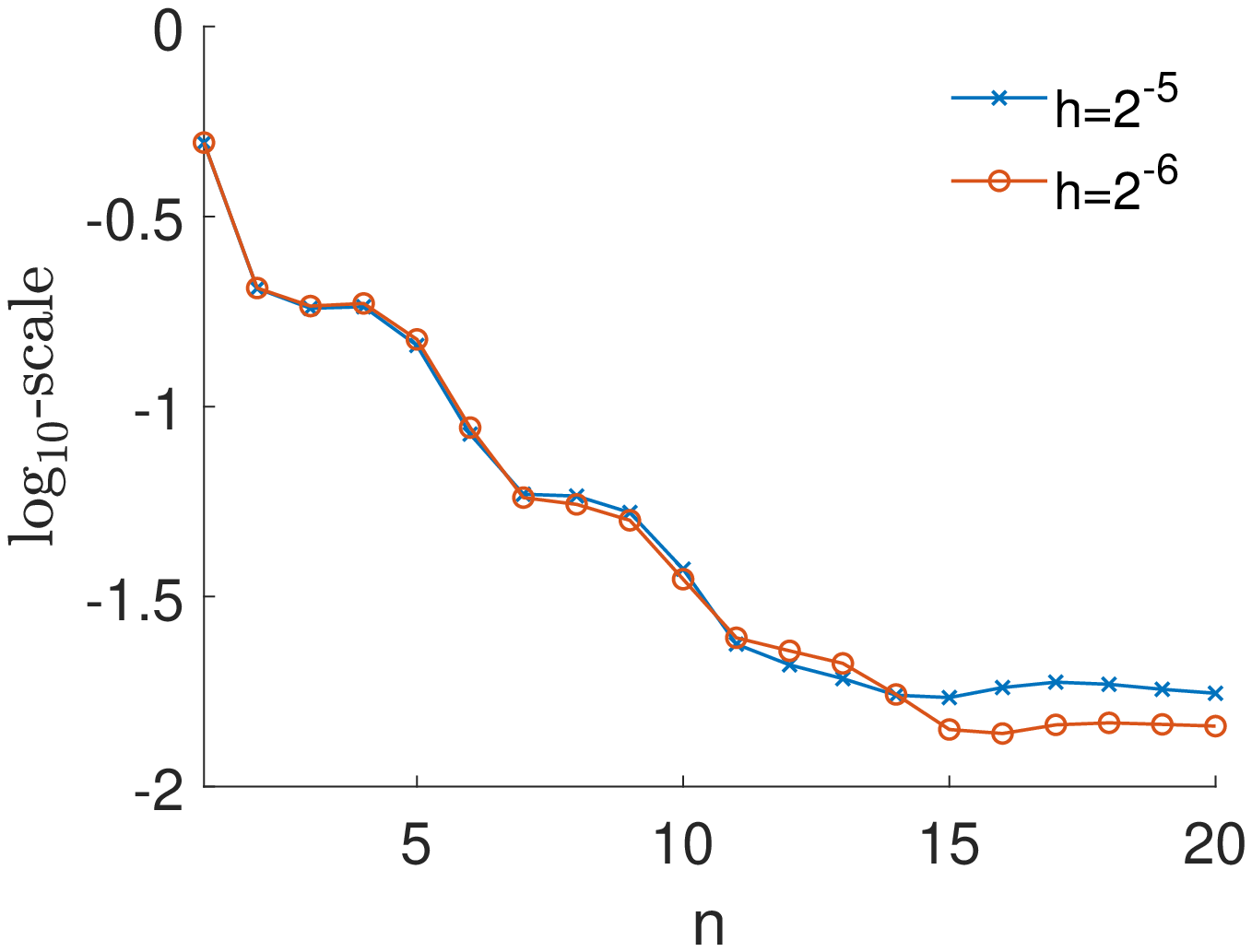}
    \end{tabular}
    \caption{Experiment 5: error behavior in terms of the iterations for the smooth load with fixed angle scattering data. In the left hand simulation $\lambda=2,\mu=1$ and $K>1$, and in the right one $\lambda=-1.1$, $\mu=1$ and $K<1$.} 
    \label{fig:ex6}
\end{figure}

Figure \ref{fig:ex6} contains the same simulation but with fixed angle scattering data both when $K\geq 1$ (left) and $K\leq 1$ (right). In both cases the load is the smooth one but Lamé parameters are different. 

\textbf{J.A. Barcel\'{o} \newline M${}^2\!$ASAI, ETSI de Caminos, Canales y Puertos,\\ 
Universidad Polit\'{e}cnica de Madrid \newline28040 Madrid, Spain \newline E-mail:
juanantonio.barcelo@upm.es }

\textbf{C. Castro  \newline M${}^2\!$ASAI, ETSI de Caminos, Canales y Puertos,\\ 
Universidad Polit\'{e}cnica de Madrid \newline28040 Madrid, Spain \newline E-mail:
carlos.castro@upm.es \bigskip}

\textbf{M.C. Vilela  \newline M${}^2\!$ASAI, ETSI de Navales,\\ 
Universidad Polit\'{e}cnica de Madrid \newline28040 Madrid, Spain \newline E-mail:
maricruz.vilela@upm.es \bigskip}
\end{document}